\documentclass[11pt,reqno,twoside]{article}

\usepackage{amsfonts,amsmath,amssymb}
\usepackage{mathrsfs,mathtools}
\usepackage{enumerate}
\usepackage{hyperref}
\usepackage{esint}
\usepackage{graphicx}
\usepackage{bm}
\usepackage{commath}
\usepackage{esint}
\DeclareMathAlphabet{\mathpzc}{OT1}{pzc}{m}{it}
\usepackage{rotating}

\usepackage{placeins} 
\usepackage{flafter} 

\hypersetup{
    hidelinks	
}

\usepackage[textsize=small]{todonotes}
\newcommand{\HA}[1]{{\color{blue}{#1}}}

\usepackage{caption}
\usepackage{subfigure}

\usepackage[absolute,overlay]{textpos}
\usetikzlibrary{positioning}

\numberwithin{equation}{section}

\usepackage[dvips,margin=1in]{geometry}

\makeatletter
\def\eqnarray{\stepcounter{equation}\let\@currentlabel=\theequation
\global\@eqnswtrue
\tabskip\@centering\let\\=\@eqncr
$$\halign to \displaywidth\bgroup\hfil\global\@eqcnt\z@
  $\displaystyle\tabskip\z@{##}$&\global\@eqcnt\@ne
  \hfil$\displaystyle{{}##{}}$\hfil
  &\global\@eqcnt\tw@ $\displaystyle{##}$\hfil
  \tabskip\@centering&\llap{##}\tabskip\z@\cr}

\def\endeqnarray{\@@eqncr\egroup
      \global\advance\c@equation\m@ne$$\global\@ignoretrue}

\usepackage{hyperref}

\usepackage[most]{tcolorbox}
\usepackage{setspace}

\usepackage{multirow}
\usepackage{colortbl}
\usepackage{tabularx}

\newcolumntype{Z}[0]{>{\centering\arraybackslash}X}%

\usepackage{authblk}

\title{Mathematical Opportunities in Digital Twins (MATH-DT) \thanks{Mathematical Opportunities in Digital Twins (MATH-DT) workshop was supported by NSF DMS - 2330895. \\
Website: \url{https://mathdt.science.gmu.edu} }
}
\date{\vspace{-5ex}}

\author{Harbir Antil} 
\affil{The Center for Mathematics and Artificial Intelligence
	and Department of Mathematical Sciences, George Mason University,
	Fairfax, VA 22030, USA.\\
    Email: hantil@gmu.edu 
    }

\begin{document}

\maketitle

\thispagestyle{empty}

\tableofcontents

\newpage
\setcounter{page}{1}

\section{Summary of Discussions}
\label{s:recommd}

The report below describes the discussions that occurred at the Workshop on Mathematical Opportunities in Digital Twins (MATH-DT) that took place at George Mason University from December 11-13, 2023.

The report highlights the `Mathematical Opportunities and Challenges'  associated with  Digital Twins (DTs). It illustrates that foundational Mathematical advances are required for DTs that are different from traditional approaches. 
The starting point of a traditional model is a generic physical law, expressed as equations, leading to a model which is typically a simplification of reality. This notion applies to problems in biology, physics, engineering or medicine. 
On the other hand, the starting point of a digital twin is a   specific ecosystem, object or person (personalized medicine, for example) which represents reality. This may require multi-scale, multi-physics modeling and coupling of various models. 
Thus, these processes begin at opposite ends of the simulation and modeling pipeline and require significantly different reliability criteria and uncertainty assessments. Furthermore, the primary goal of a DT is to assist humans in making decisions for the physical system, which (via sensors) in turn feeds data into the DT, and the cycle continues throughout the life of the physical system. Again, this is unlike the existing traditional approaches.

While some of the foundational mathematical research can be done without a specific application context, one must also keep specific applications in mind for DT research. For instance, modeling a bridge is very different from modeling a biological system, such as a patient, or a socio-technical system, such as a city. The models are different, ranging from differential equations (deterministic or uncertain) in most engineering applications, to various stochastic models in biology, including agent-based models. They will most commonly be multi-scale hybrid models or large scale (multi-objective) optimization problems under uncertainty 
(e.g. parameter identification). There are no universal models or approaches. Kalman filter methods for forecasting, for instance, might work well in many engineering applications, but tend to do very poorly in biomedical applications, even for ODE systems, and have not even been seriously explored for agent-based models. Ad hoc studies have shown that Artificial Intelligence/Machine Learning (AI/ML) methods can fail for simple engineering systems and can work very well for biomedical applications, but only limited systematic work has been done. In summary, the DT research presents many mathematical challenges and opportunities as highlighted throughout the report and summarized in  section~\ref{s:recom}. 

\vspace{-0.5em}
\begin{itemize}\itemsep0pt

\item {\bf Direction 1.} Clarify in the mindset of academia, computational mechanics software companies, industry, funding agencies, and potential users that foundational Mathematical advances are required for DTs that are different from traditional fields. In order to make DTs a reality, considerable advances will be required, for exmaple, in large-scale optimization, inverse problems, data assimilation, multi-scale modeling, multi-physics modeling, uncertainty quantification, model order reduction, scientific machine learning, quality assessment, process control, rapid statistical inference, online change detection for actionable insights, software development (see section~\ref{s:chal} and \ref{s:soft}), all in combination and in the context 
of particular DTs. Some of the methods may only work well in certain fields - some may be general. 

    \item {\bf Direction 2.} It is crucial to develop benchmark test problems and data sets for different applications. Examples are: bridges (in civil engineering), vibrating beams (in mechanical engineering), and sepsis (in medicine). See section~\ref{s:opp} for more examples. It is also critical to develop mathematical and statistical principles for DT model updating and maintenance.

    \item {\bf Direction 3.} Mathematics research for digital twin applications needs to be done in an interdisciplinary context. Research funding announcements (RFAs) 
    need to reflect that. In many cases, the appropriate approach is collaborative RFAs between different divisions, directorates, and agencies. Funding agencies should organize team-building events that connect mathematicians with domain experts. The NSF has some experience with this (\emph{Ideas-Lab}, for instance), and so does NIH and DOE. 

    \item {\bf Direction 4.} NSF DMS is recommended to fund a unique DT collaborative institute in the mold of the mathematics institutes, that can provide a collaboration platform, and a repository of resources. 
    Notice that MATH-DT workshop was different than a traditional workshop and it observed a significant participation from government funding agencies
    (e.g., NSF, NIH, DOD, and DOE), industries (e.g., Siemens, NVIDIA, ESI), and national labs (e.g., Sandia National Lab, Lawerence Livermore National Lab, Argonne National Lab, US Naval Research Lab, Brookhaven National Lab). The workshop has set a stage to create this  collaborative mathematics institute between Academia, National Labs, and Industry. Notice that nothing like this currently exists.

    \item {\bf Direction 5.} It is recommended to create RFAs for research institutes, similar to the AI institutes, as well as various consortia focused on DT technology. Significant focus should be on interdisciplinary aspects, critical to DT success, consisting of mathematicians, statisticians, data scientists and domain experts.  There are several European models that may also serve as templates. 

	\item {\bf Direction 6.} There should be RFAs for training programs, summer schools, tutorials, and development of curricular materials. 
 For instance it is strongly recommended to create an annual (1-2 weeks long) summer school on this topic. During the summer school, the students could focus on a few specific benchmark problems in DTs and also develop abstract mathematical skills around those problems. Indeed, these schools could be a excellent opportunity to bring different worlds together which are essential to DTs, e.g., theory and practice, data analysis and data creation, simulation and experiment.
 One expected outcome could be to produce a Volume every year, which is organized in the form of a textbook. 
		An example of this exists already \cite{HAntil_DPKouri_MDLacasse_DRidzal_2018a}

  Another discussion was to create RFAs for interdisciplinary fellowships. An example of this is Department of Energy Computational Science Graduate Fellowship (DOE CSGF). 
  \end{itemize}

We acknowledge two additional efforts on Digital Twins. In particular, the workshop on ``Crosscutting Research Needs for Digital Twins" at the Santa Fe Institute (October 11-12, 2023) \cite{santafe} and the National Academy of Sciences report on ``Foundational Research Gaps and Future Directions for Digital Twins" \cite{national2023foundational}. The current report is complementary to these other two reports, with emphasis on Mathematics Opportunities in Digital Twins. We further emphasize that some of our suggestions mirror the conclusions and recommendations of the other two reports.

\newpage

\section{Overview and Goals}

\subsection{Background}
\label{s:back}

Mathematics is undoubtedly the language of sciences such as Physics, Chemistry, 
Mathematical Biology, and Data Science. 
Mathematics is more important than ever, especially with the emergence of Artificial 
Intelligence and Machine Learning (AI/ML) as one of the priority research areas. 
In comparison to Mathematicians, more Computational Engineers are working in these
areas and are regularly using fundamental concepts learned from decades of investment in 
mathematics, in particular in optimization, numerical analysis, partial differential 
equations, numerical linear algebra, statistics and scientific computing. 

As it happened with other research directions, such as finite element analysis, one 
cannot bring a significant change, or explain the AI/ML concepts in the required depth 
or prove convergence of widely used algorithms, without new mathematical developments.
Absence of mathematics can also impede the development of new algorithms. Especially 
for AI/ML, the situation is alarming. For instance, purely 
data driven AI/ML approaches can fail to generalize on basic physics based models,
problems that can be solved within milliseconds using traditional Finite Element 
based approaches. Additionally, most of the fundamental questions remain open, 
for instance, what kind of network should be used, the size (depth and width) of the network, how much training data is sufficient, how to train highly nonconvex 
optimization problems, what to expect in terms of accuracy and when 
one can trust these approaches. In many cases it is not even clear how to rigorously 
phrase these questions, especially for physics problems `in the absence of physics'. 
So clearly, mathematics has a significant role to play here.  

Let us now comment on what led to the emergence of AI/ML, especially in science and engineering.
Indeed, this hasn't happened in a vacuum. Traditionally, the focus has been on model based simulations 
which have their own limitations. It was realized that not everything can be modeled using first 
principles. Additionally, there were significant limitations in terms of computing resources. The 
emergence of advanced computing resources has fast-tracked developments in AI/ML to address these limitations. Many times this came at the expense of fundamental physical principles. 

Now we may ask whether one could reap the benefits of all these advancements. Additionally, notice
that the starting point of a traditional model is a generic
physical law, expressed as equations, leading to a model which is typically a simplification of reality.
This notion applies to problems in biology, physics, engineering or medicine. A different thinking is needed if one starts from a specific ecosystem, 
object or person in a more realistic setting. Indeed, these processes begin at the opposite ends of a pipeline and require significantly different reliability criteria 
and uncertainty assessments. 
In this vein, the AIAA 2020 Institute Position Paper \cite{aiaa2020digital} defines the term `Digital Twin (DT)' as  
\begin{center}
``\emph{A Digital Twin is a set of virtual information constructs that mimics the 
structure, context, and behavior of an individual/unique physical asset, is dynamically 
updated with data from its physical twin throughout its lifecycle, and informs decisions 
that realize value.}"
\end{center}
This definition is leading to a fundamental shift on how we envision mathematics and science. 

Why have DTs emerged recently? As is the case so often in disruptive changes, because of several megatrends: a) the pervasive use of Computer Aided Design (CAD) systems, b) the widespread availability and use of computational tools to `pre-compute, only then build' and `pre-compute, then operate'; and c) the emergence of precise, rugged, connected and cheap sensors and cameras that may be used to monitor products, processes and patients.

Indeed, DT is truly a completely new interdisciplinary framework which encompasses traditional
simulation based approaches, data science, optimization, inverse problems, uncertainty 
quantification, and AI/ML with applications to a wide range of problems. It requires researchers of 
various backgrounds (e.g., domain scientists, computational engineers, etc.) working together 
in both small and large teams with mathematicians {\bf both at the core and as a bridge}. 

\begin{figure}[h!]
	\centering	
	\includegraphics[width=\textwidth]{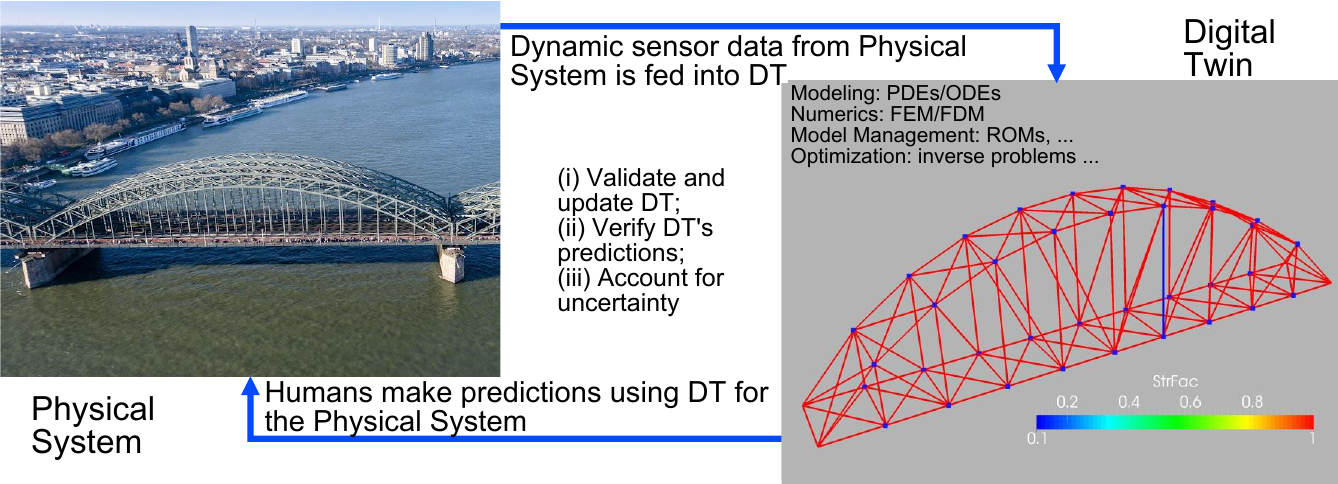}
	\caption{\label{f:bridge}
	\small 
	An actual bridge (physical system) is shown on the left and the digital twin (DT)
        of this bridge is shown on 
	the right. From the sensors on the real bridge, a dynamic dataset is collected. The relevant 
	dataset is identified (possibly using AI/ML techniques) and fed into the DT. The DT 
	involves multiple components, depending on the purpose. The goal here 
	is to identify a weak beam marked in blue color while accounting for the sensor data 
	and the underlying physics (elasticity equations). Then the DT informs an engineer to 
	make decisions about the actual bridge to make appropriate fixes. These decisions 
	should account for uncertainty and must be risk-averse. 
	This entire problem can be cast as a minimization problem with tracking the displacement or
	strain measurements subject to elasticity equations as constraints. To create a risk-averse
	framework, one could additionally consider uncertainty in loads and minimize risk measures,
	such as Conditional Value at Risk (CVaR). The impact of CVaR is illustrated in 
        section~\ref{s:bridge} on a crane.}
\end{figure}

A non-exhaustive list of examples include: DTs of bridges and high-rise buildings, 
DTs of nuclear reactors, conventional power stations and renewable energy parks, 
DTs of production and manufacturing processes, DTs of humans (e.g., for biomedical interventions or shopping preferences). 
The primary goal of a DT is to assist humans to make decisions for the physical system, 
which in turn feeds data into the DT and the cycle continues. Figure~\ref{f:bridge} shows an actual 
bridge. Notice that, the DT framework is unique and none of the existing frameworks consider such a setup; and a setup for the entire lifetime of a bridge. 
Additionally, as pointed out earlier, the starting point for a digital twin is a specific ecosystem, 
object or person in a more realistic setting. This is unlike a model, which is a simplification of reality, where the starting point is a generic physical law often expressed as equations. Thus, the processes in DTs begin at opposite ends of a 
pipeline and require significantly different reliability criteria and uncertainty assessments. 

The mathematical opportunities and challenges in DTs are manifold with a significant broader impact. 
See section~\ref{s:examples} for more details and examples. 

\subsection{MATH-DT Workshop}

The workshop on Mathematical Opportunities in Digital Twins (MATH-DT) took place  
at George Mason University (GMU), Arlington Campus from December 11-13, 2023 to 
determine how Mathematics can contribute to DT research and how DT can help Math research. 
This workshop was primarily supported by NSF DMS 2330895. The workshop also received 
a partial support from the SIAM Washington-Baltimore Section.

\smallskip
\noindent
\emph{Organizers:}
MATH-DT was chaired by Harbir Antil (GMU). The other co-organizers were Benjamin 
Seibold (Temple University) and Kathrin Smetana (Stevens Institute of Technology). 
Locally, it was also supported by Ratna Khatri (U.S. Naval Research Laboratory),  
Rainald L\"ohner (GMU), and Roland W\"uchner (Technical University Braunschweig, 
Germany). 

\begin{figure}[htb!]
\centering
\includegraphics[width=0.6\textwidth]{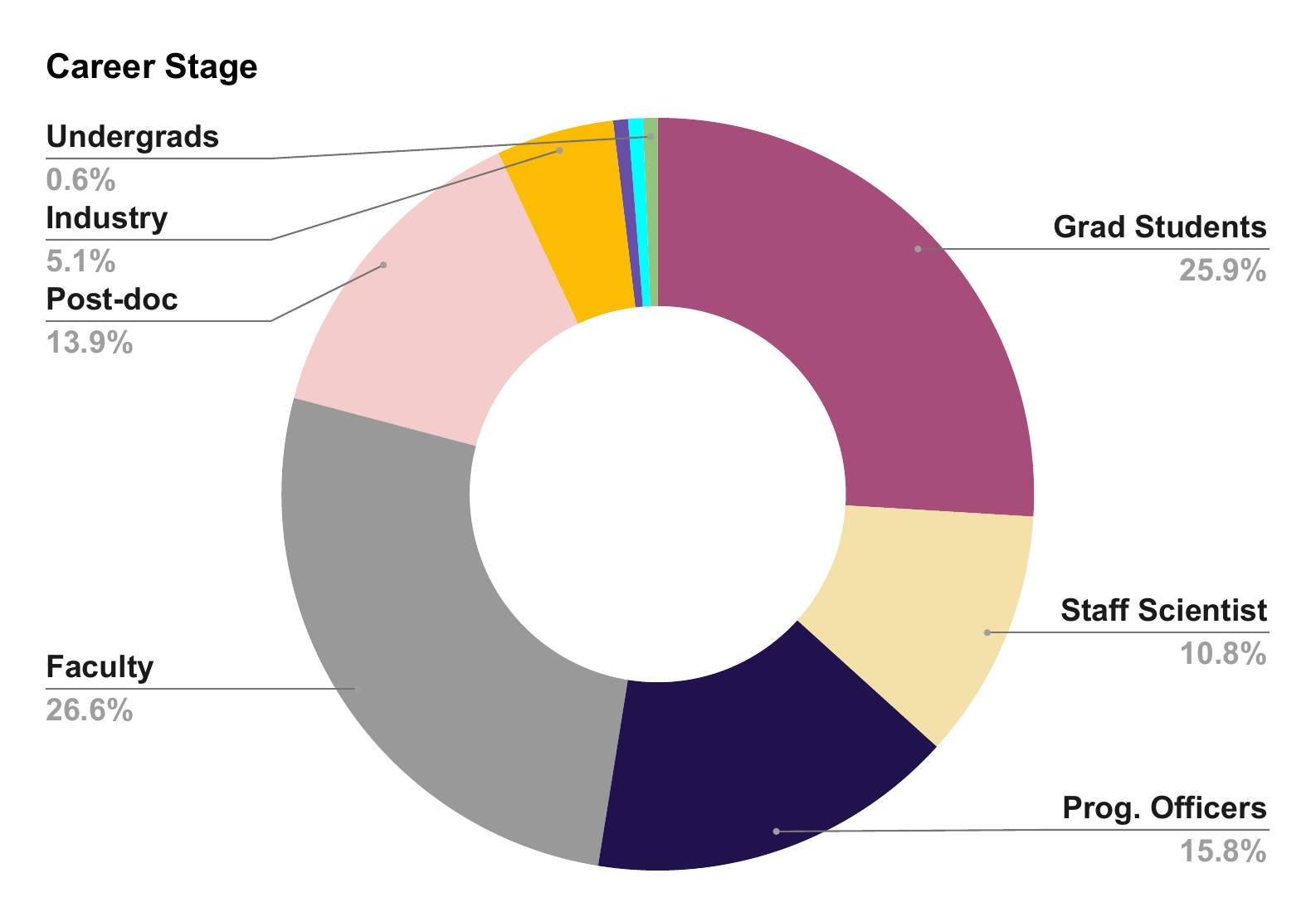}
\caption{\label{f:dist}The panel shows the distribution of MATH-DT attendees.}
\end{figure}

\begin{figure}[h!]
    \centering
    \includegraphics[width=0.8\textwidth]{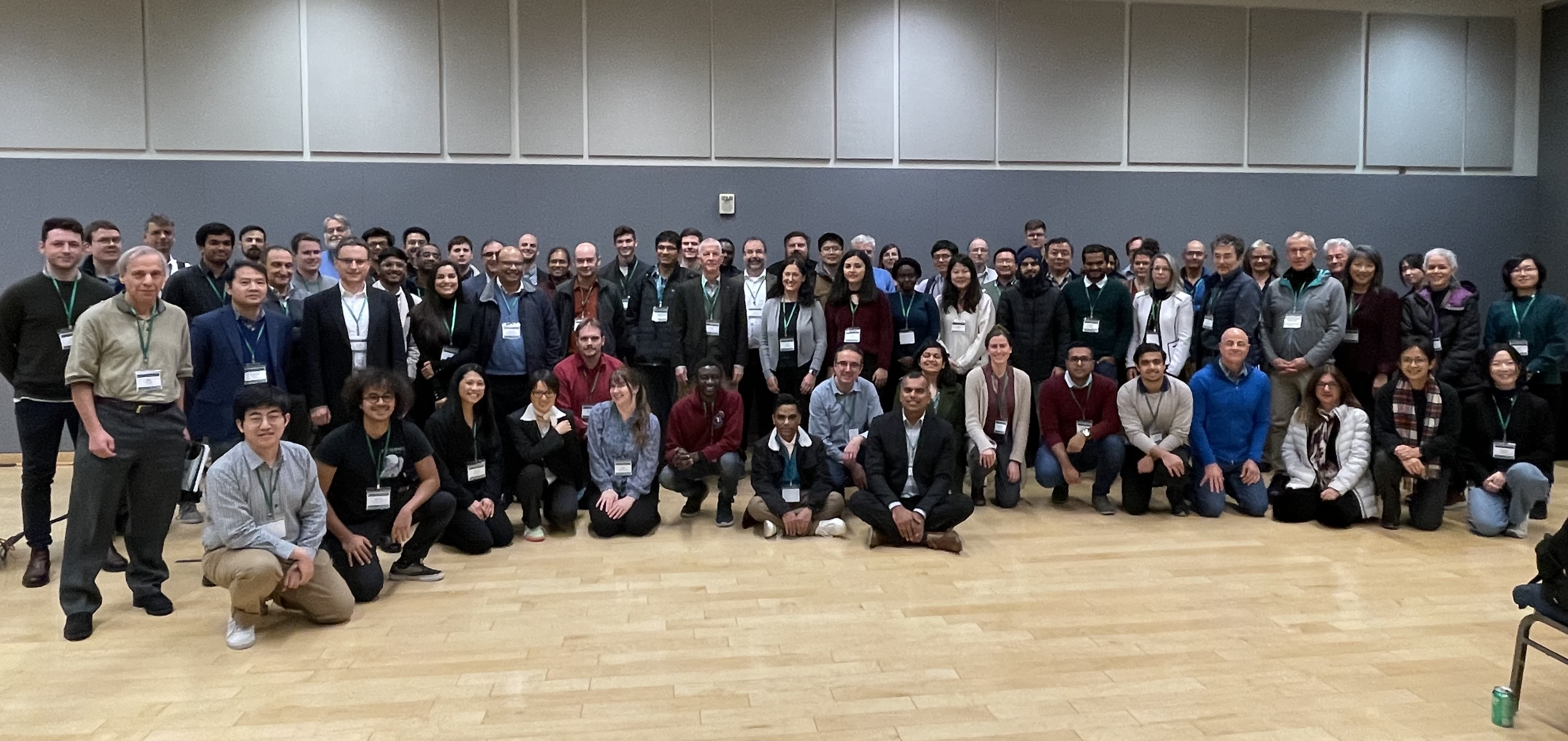}    
    \caption{Some of the attendees of MATH-DT Workshop}
    \label{f:attendees}
\end{figure}

\smallskip
\noindent
\emph{Participants:}
The workshop received 180 registrations and it was attended by 128 people in-person 
with additional 15 people online over zoom. Figure~\ref{f:dist} shows a distribution 
of registered participants. Notice a large representation from students and postdocs 
(around 40\%). 
A significant participation was also observed from government funding agencies
(e.g., NSF, NIH, DOD, and DOE), industries (e.g., Siemens, NVIDIA, ESI), and 
national labs (e.g., Sandia National Lab, Lawerence Livermore National Lab, Argonne 
National Lab, US Naval Research Lab, Brookhaven National Lab). 
Some of the workshop attendees are also shown in Figure~\ref{f:attendees}.
The full workshop 
schedule can be found in Appendix~\ref{s:schedule}, the list of speakers and panelists 
is available in Appendix~\ref{s:speakers}. All the titles and abstracts are available 
under Appendix~\ref{s:abstracts}.

 \subsection{Goals and Format}
 
The primary goal of the workshop was to identify how Math can
benefit a DT and how a DT can benefit Math. In particular, MATH-DT was 
interested in posing `Hilbert's 23' style problems with emphasis on:
\begin{enumerate}[(Q1)]\itemsep0pt
	\item What can Math do for Digital Twins (DT) ({\bf mathematical opportunities})?
		\begin{enumerate}[(a)]
			\item What is the mathematical state-of-the-art on Digital Twins in your field ?
			\item What impact can you foresee for an increased use and rigor of mathematics on Digital Twins ?
		\end{enumerate}
	
	\item What can DT do for MATH ({\bf mathematical challenges}) ?
		\begin{enumerate}[(a)]
			\item What are the current mathematical limitations/challenges ?
			\item What mathematical tools and developments are needed to overcome these limitations/challenges ?
		\end{enumerate}
		
	\item Broader impact
		\begin{enumerate}[(a)]
 			\item What impact will Digital Twins have on your field ?
			\item What impact will Digital Twins have on society at large ? 
		 \end{enumerate}	
\end{enumerate}

\medskip
\noindent
\emph{Format:} The workshop consisted of 12 plenary talks, 6 application focused talks,
one funding agency panel, and 30 poster presentations from early career researchers 
(graduate students and postdocs). In addition, each day there was a `Dream Big (DB)' session
where all the participants took part. Five topics were identified under this session. All 
the speakers were asked to address (Q1)-(Q3) within the context of these five topics,
see Section~\ref{s:dream} and Appendix~\ref{app:DB} for the group findings. The overall
summary is provided in section~\ref{s:recom}. 
\begin{enumerate}[(DB1)]\itemsep0pt
	\item Data
	\item Modeling and Forward Problems 
	\item Optimization and Inverse Problems 
	\item Validation, Verification and Predict
	\item Software 
\end{enumerate}
%

\section{Examples of Digital Twins and Their Mathematical Formulation}
\label{s:examples}

In this section, several examples of DTs are provided. The fist example corresponds to neuromorphic imaging, where real data from the International Space Station (ISS) \footnote{We acknowledge the neuromorphic data provided by the Space Physics and Atmospheric Research Center (SPARC) at the United States Air Force Academy.  SPARC built and flew the Falcon Neuro mission to the International Space Station in collaboration with Western Sydney University (WSU).} has been used. In the second example, location of cracks / weakness in structures is studied while accounting for uncertainty in loads. The third example corresponds to medical digital twins.

\subsection{Neuromorphic Imaging}
Neuromorphic imaging, embodied by the evolution of Event Cameras or Dynamic Vision 
Sensors (DVS), introduces a transformative approach to visual information acquisition. 
Departing from traditional frame-based cameras, neuromorphic cameras capture visual 
data asynchronously, responding exclusively to intensity changes at each pixel location. 
Unlike fixed-rate sampling in traditional cameras, neuromorphic cameras adapt their 
sampling rate based on scene dynamics, avoiding under-sampling of swiftly changing 
scenes or redundant over-sampling of slowly changing ones.

They present an interesting avenue for DTs. Consider an example where the goal is to 
capture dynamic images of earth from the fast moving International Space Station
(ISS). This is an extremely high contrast environment with limited power. 
DT is constructed using neuromorphic (event) data as illustrated in Figure~\ref{f:neuro} by solving an inverse 
problem of type \cite{HAntil_DSayre_2023a} 
\begin{equation*}
\bm{v}_{xy}=\operatorname*{argmin}_{\bm v_{xy}} \frac{1}{2} \| {\bm A} \bm v_{xy}  -  \bm E_{xy} \|_2^2 + \frac{\lambda}{2} \| \bm v_{xy}\|_2^2 \, .
\end{equation*}
Notice that the feedback from the DT can, for instance, help adjust the camera position.
		\begin{center}
			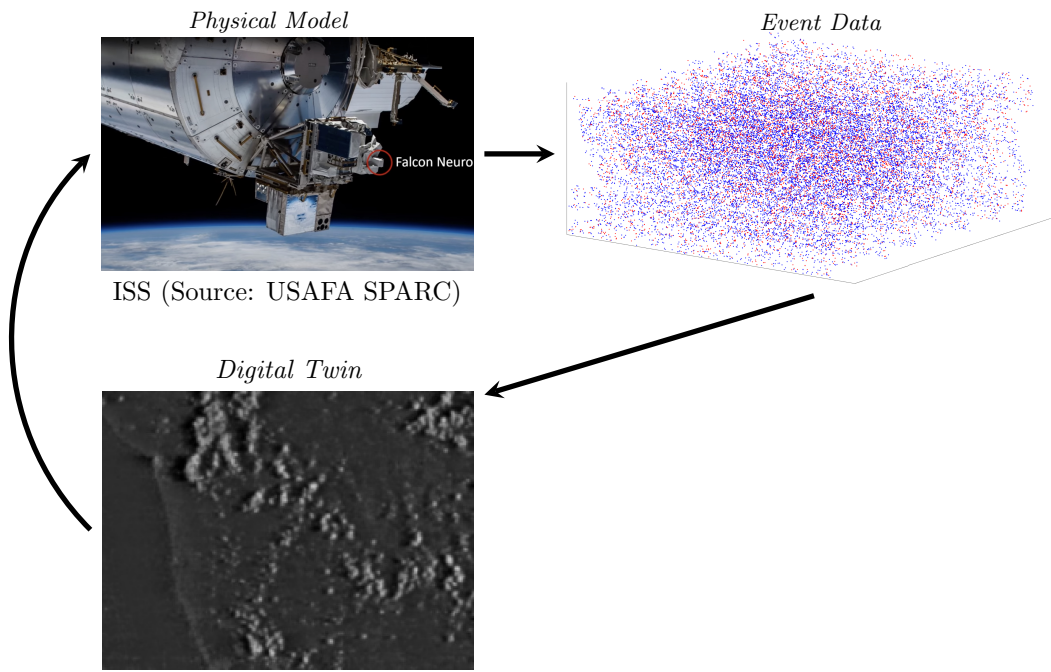
\begin{figure}[h!]
				\begin{tikzpicture}[>=stealth,thick]
				\node (img1) at (-1.0,5) {\includegraphics[width=0.3\textwidth]{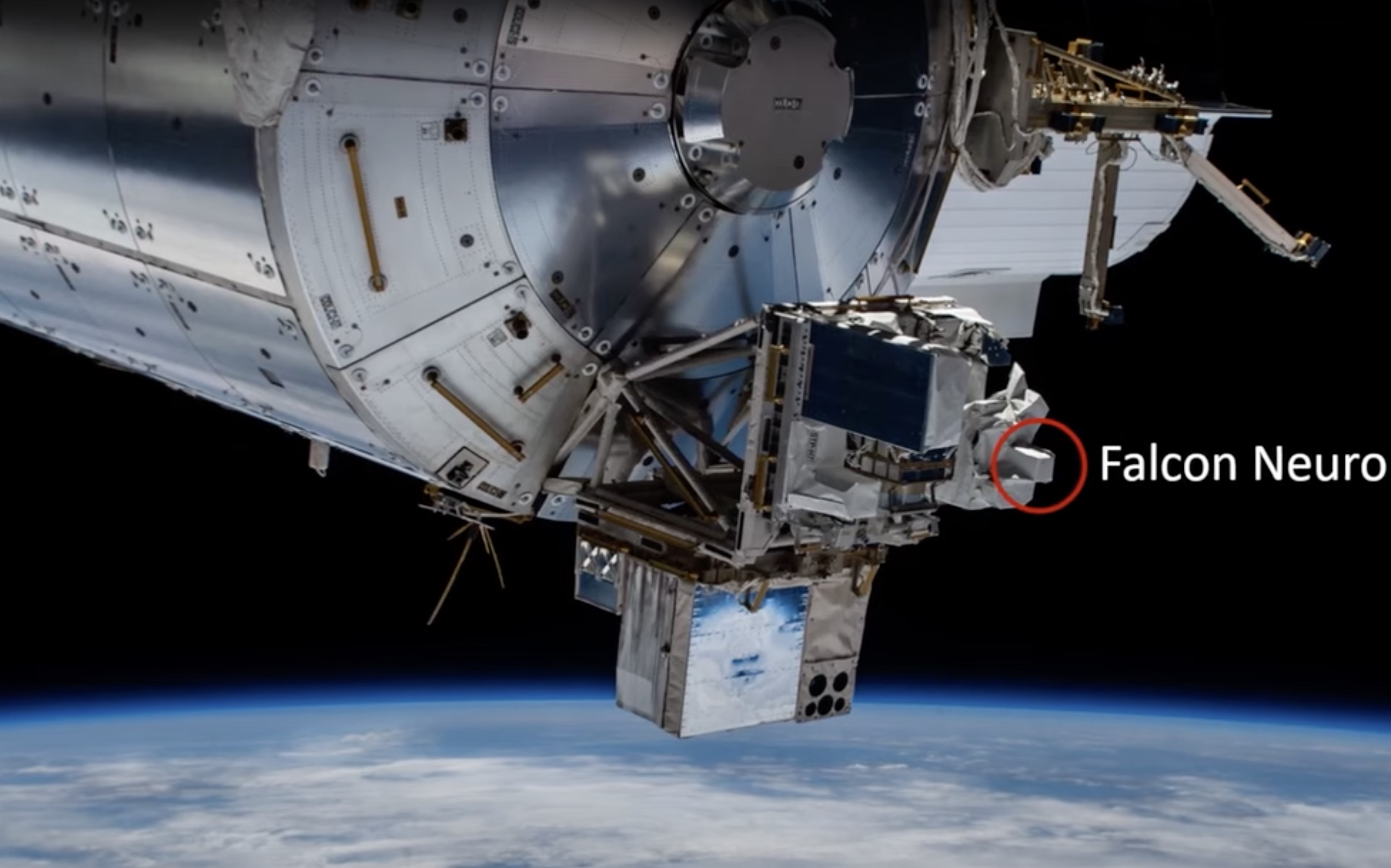}};
				\node[below=of img1, yshift=3em] {{\small ISS (Source: USAFA SPARC)}};
				\node (img2) at (6,5) {\includegraphics[width=0.4\textwidth]{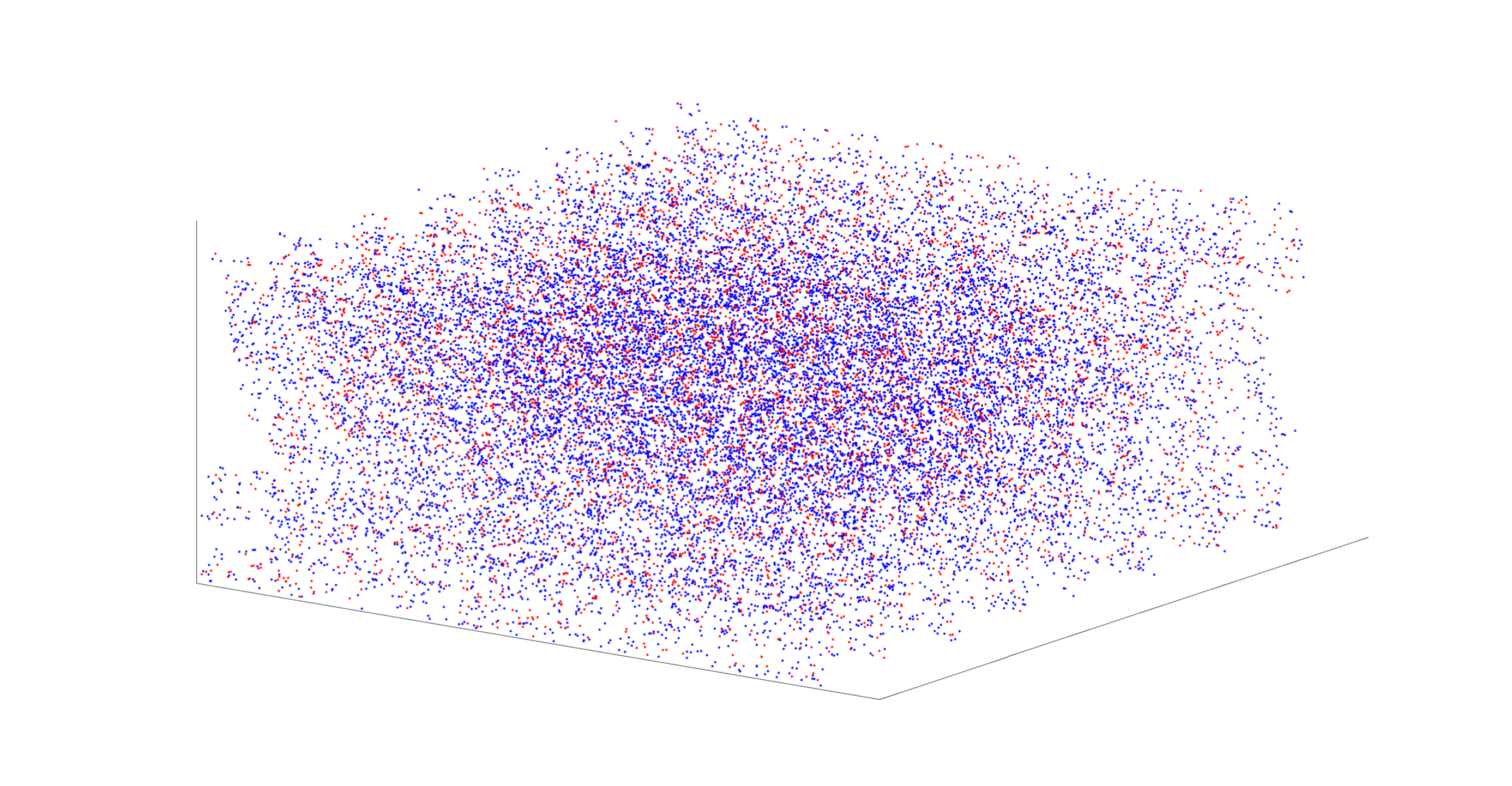}};
				\node (img3) at (-1.0,0.0) {\includegraphics[width=0.30\textwidth]{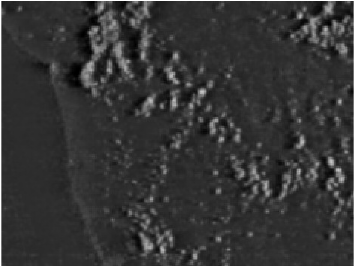}};
				
\draw[->, line width=2pt] (img1.east) -- (img2.west) node[midway, right, xshift=-130pt, yshift=50pt] {\footnotesize \emph{Physical Model}};
\draw[->, line width=2pt] (img2.south) -- ([yshift=-6pt]img3.north east) node[midway, below, xshift=65pt, yshift=130pt] {\footnotesize \emph{Event Data}};
\draw[->, line width=2pt] (img3.west) to[out=135, in=-135]  (img1.west) node[midway, left, xshift=5pt, yshift=60pt] {\small \emph{Digital Twin}};
				\end{tikzpicture}				
			\caption{\label{f:neuro} Figure shows a Digital Twin loop. In the top panel (left), a neuromorphic camera located on the international space station (ISS) 
			captures event data over a time period. This data is shown in panel on the right. One then solves an 
			inverse problem to obtain the image reconstruction. The latter helps, for instance, 
			adjust the camera position to track objects of interest more accurately and this process continues.}
			\end{figure}
		\end{center}

\subsection{Identifying Weakness in Structures: Bridges and Cranes}
\label{s:bridge}

Mathematically, the DT problem described in Figure~\ref{f:bridge} is an optimization 
problem with elasticity equations as PDE constraints. Let the physical domain be 
given by $\Omega$ and let the boundary $\partial\Omega$ be partitioned into the 
Dirichlet ($\Gamma_D$) and Neumann ($\Gamma_N$) parts, respectively. Then the 
optimization problem, to identify the weakness in the structure, i.e., coefficient function 
$z$, is given by 
		\begin{equation}\label{eq:min}
			\min_{ ( \bm u, z) \in U \times Z_{\rm ad}} 
			J(\bm u, z)  
			\quad \mbox{subject to} \quad
			-\mbox{div} (\bm \sigma(\bm u;z) ) = \bm f ,
		\end{equation}		
with $\bm u = \bm 0 \mbox{ on } \Gamma_D$ and 
$\bm \sigma(\bm u;z) \bm \nu = \bm 0  \mbox{ on } \Gamma_N$. Here 
$\bm \sigma$ is the stress tensor, $\bm f$ is the given load, 
$Z_{\rm ad} := \{ z \in L^2(\Omega) \, : \, 0 < z_a \le z \le z_b \mbox{ a.e. }
x \in \Omega \}$ is a closed convex set with $z_a, z_b \in L^\infty(\Omega)$. 
Here $J$ indicates the misfit between the displacement $\bm u$ or
strain tensor, and their respective measurements (displacement or strain). 

\begin{figure}[!hbt]
    \centering
    \subfigure[Target Displacements and Sensor Locations]{\includegraphics[width=0.45\linewidth]{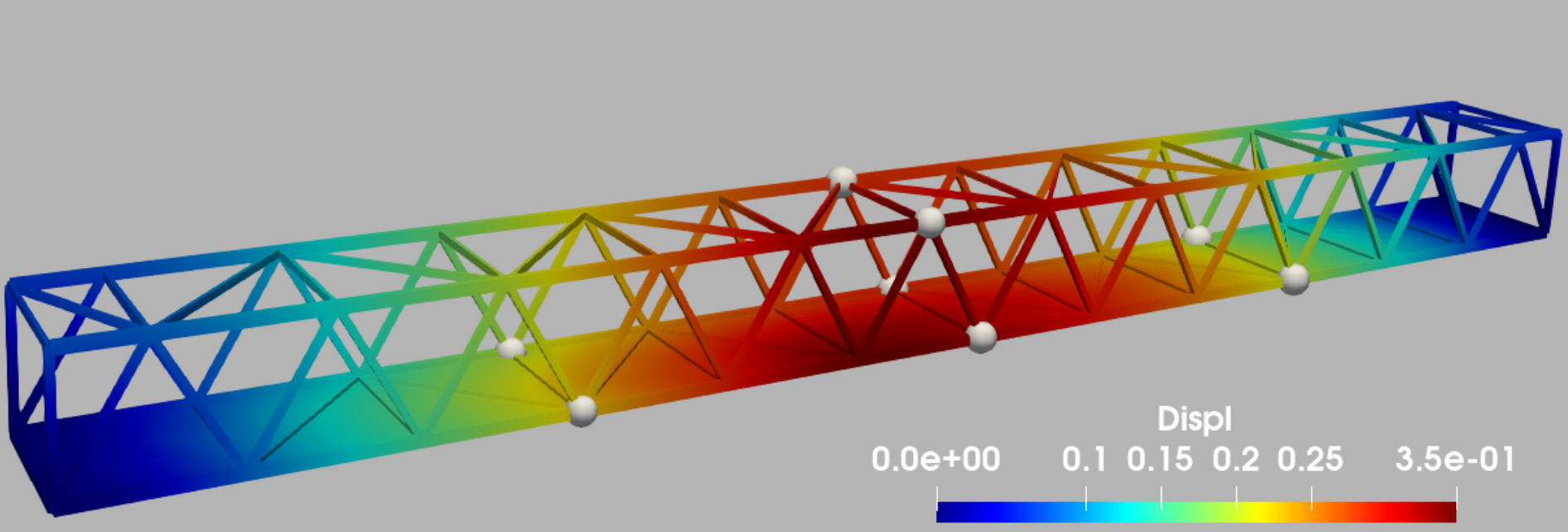}}
    \subfigure[Target Strength Factor]{\includegraphics[width=0.45\linewidth]{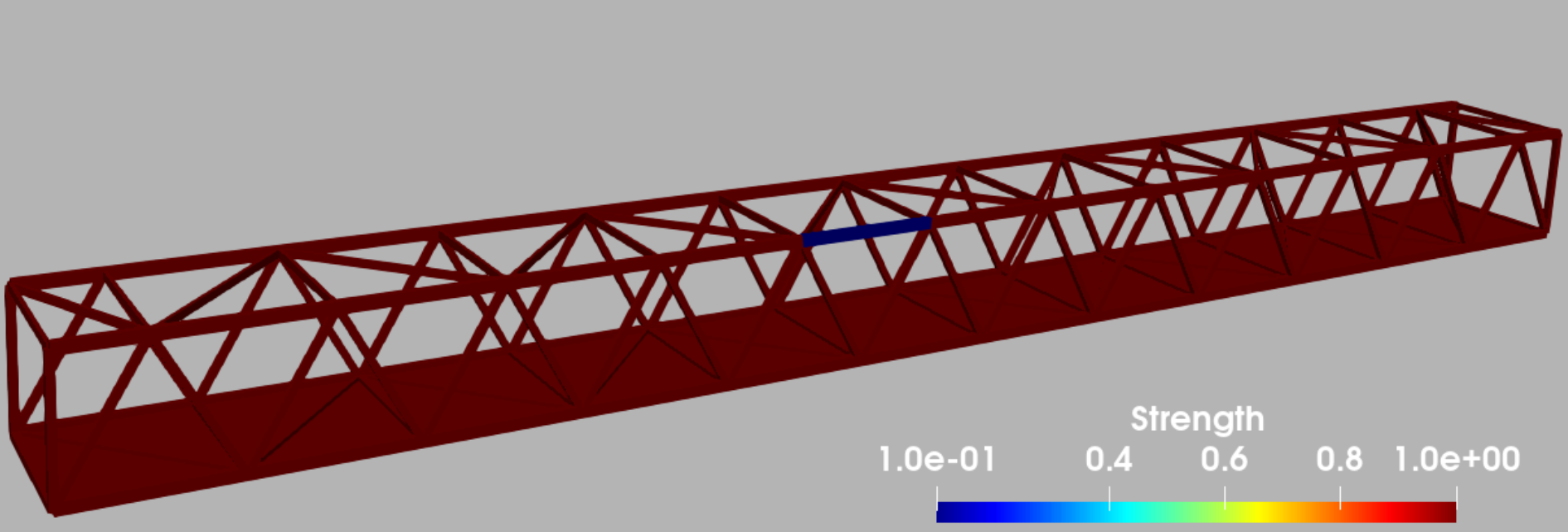}}
    \subfigure[Displacements Obtained]{\includegraphics[width=0.45\linewidth]{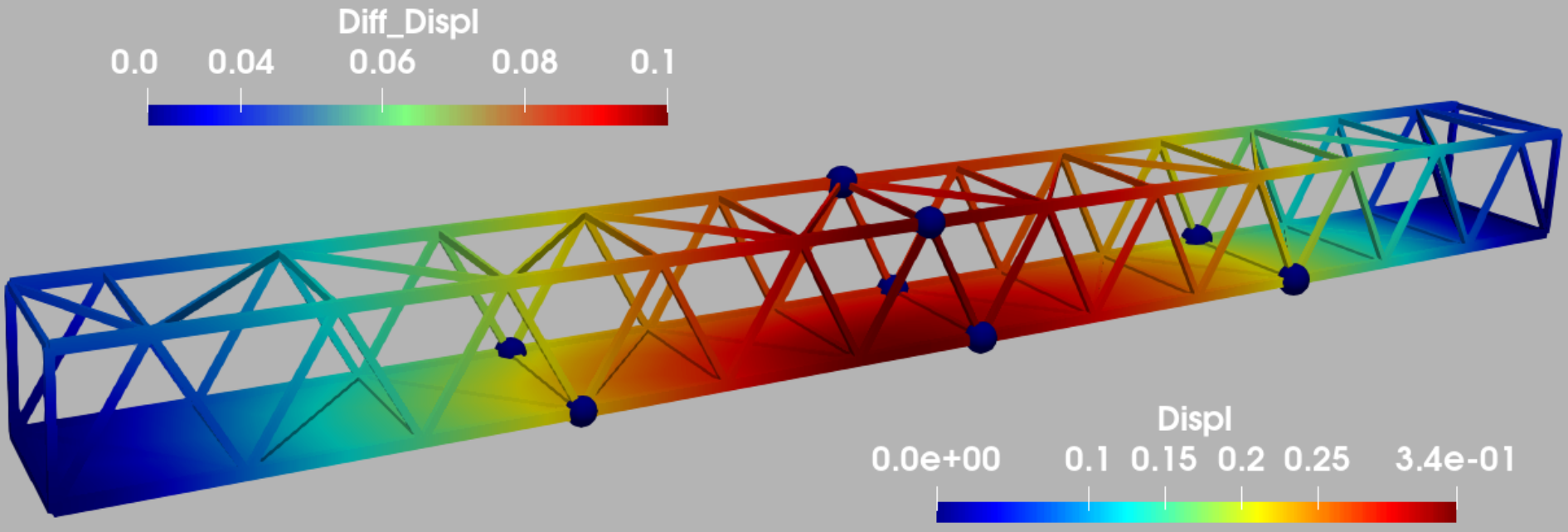}}
    \subfigure[Strength Factor Obtained]{\includegraphics[width=0.45\linewidth]{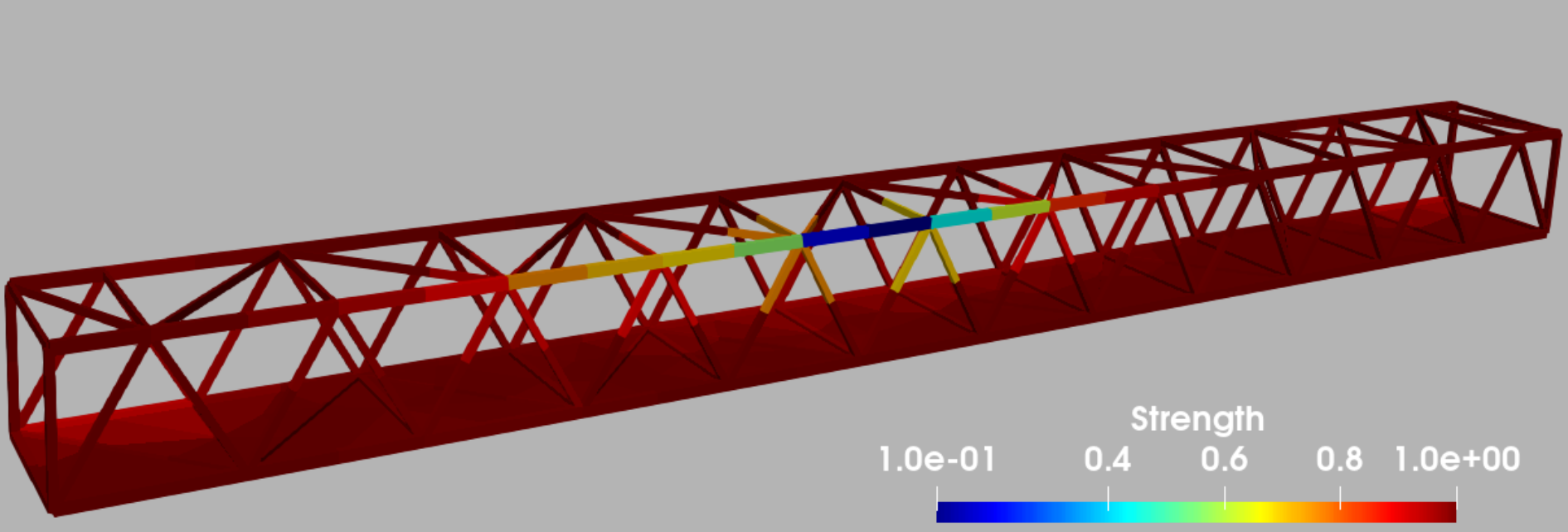}}
    \caption{Top row: Target displacement $\bm u$ (left) and strength factor $z$ (right). 
             The color bars at the bottom right corners, in each panel, correspond to the actual displacements and strength factors, 
             respectively. 
             Bottom row: Solution obtained at the 200-th optimization iteration. Panels (c) and (d), 
             respectively display optimization results. The color bar on the top left in panel (c) displays the magnitude 
             of the difference between target and actual displacements at the measuring points 
             (in \emph{m}). 
             The DT is able to identify weakness in the structure.}
    \label{fig:footbridge_target}
\end{figure}

A result is shown in Figures~\ref{fig:footbridge_target} (cf.~\cite{FAiraudo_RLoehner_HAntil_2023a}).
Panel (a) shows the target displacement and sensor locations and (b) shows the 
target strength factor $(z)$. A solution to the optimization problem \eqref{eq:min}
under various loading conditions is shown in (c) and (d) with
respective displacement and strength factor. 

\begin{figure}[!h]
    \centering
    \subfigure[Target strength]{\includegraphics[width=0.45\linewidth]{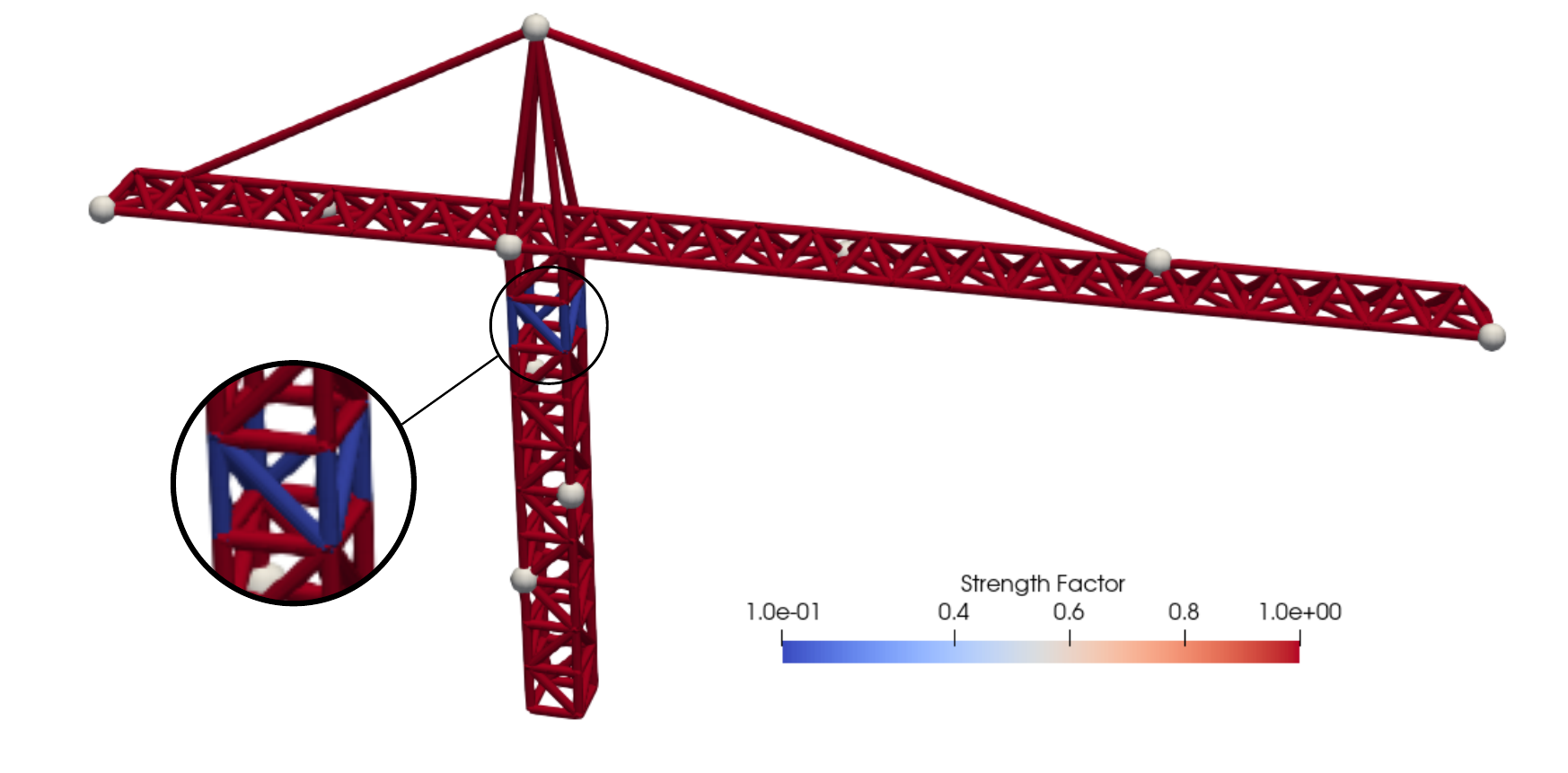}} \quad
    \subfigure[$\mathbb{E}(X)$]{\includegraphics[width=0.45\linewidth]{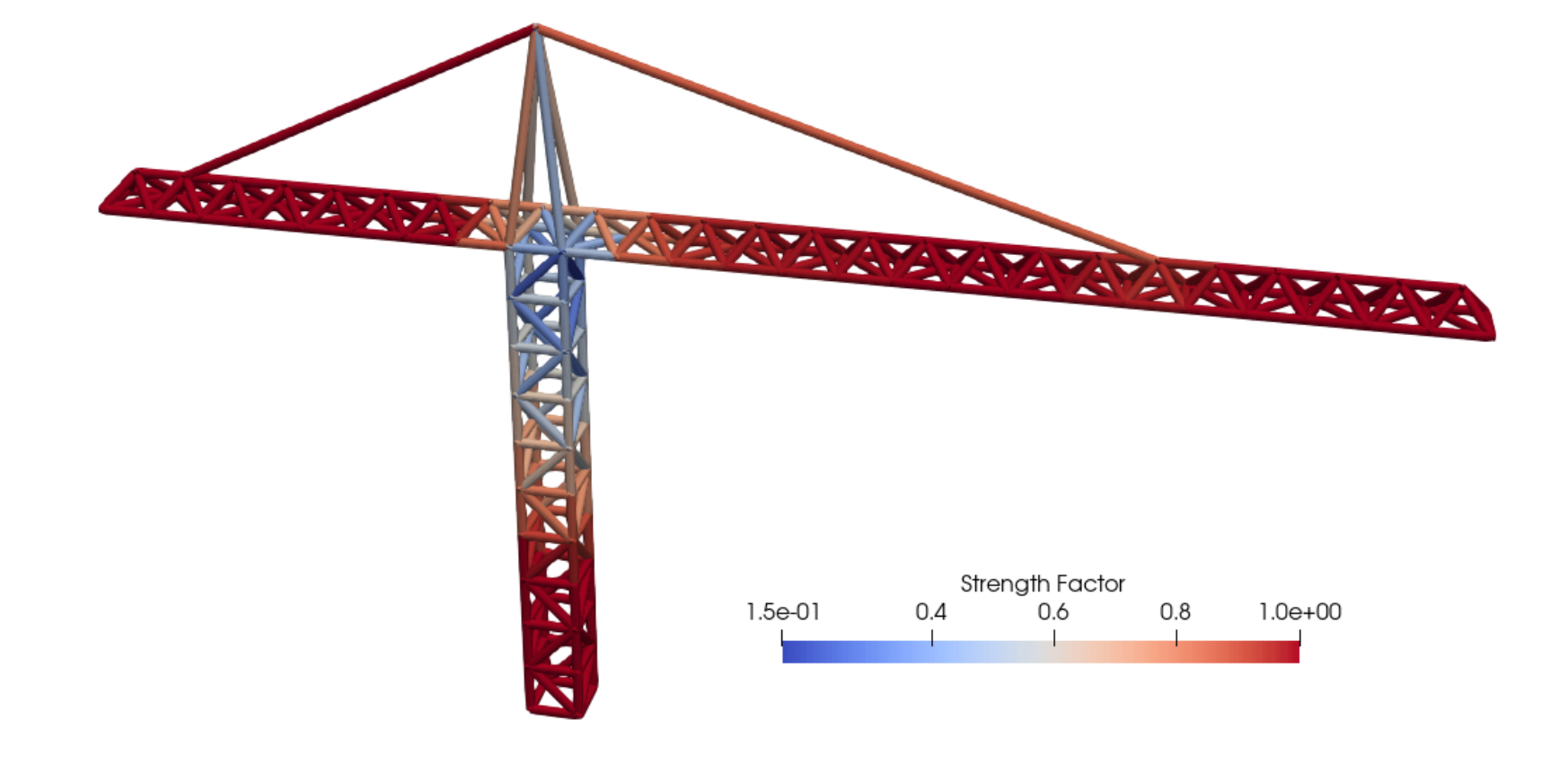}} 
    \subfigure[CVaR$_{0.3}(X)$]{\includegraphics[width=0.45\linewidth]{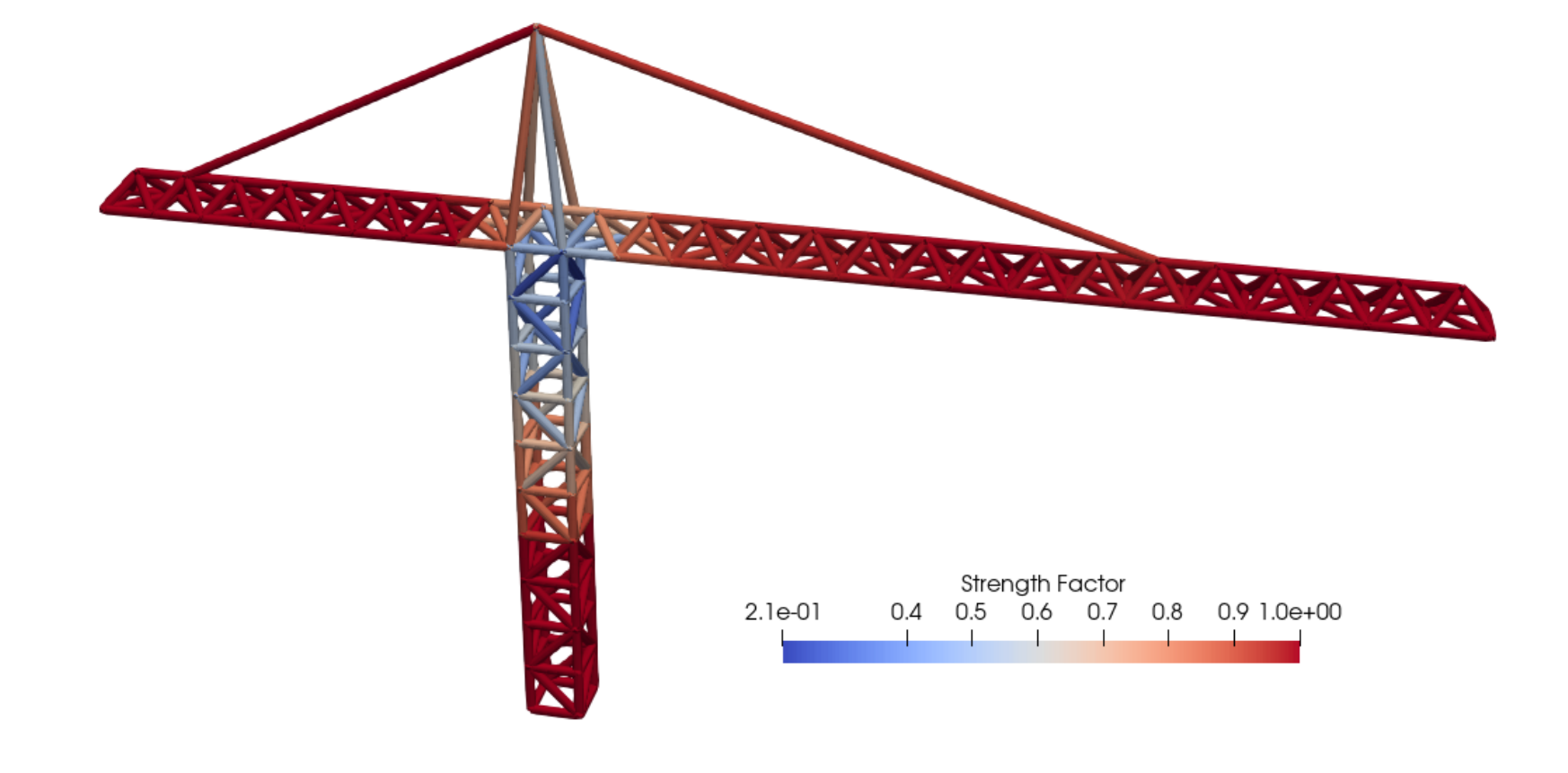}} \quad
    \subfigure[CVaR$_{0.8}(X)$]{\includegraphics[width=0.45\linewidth]{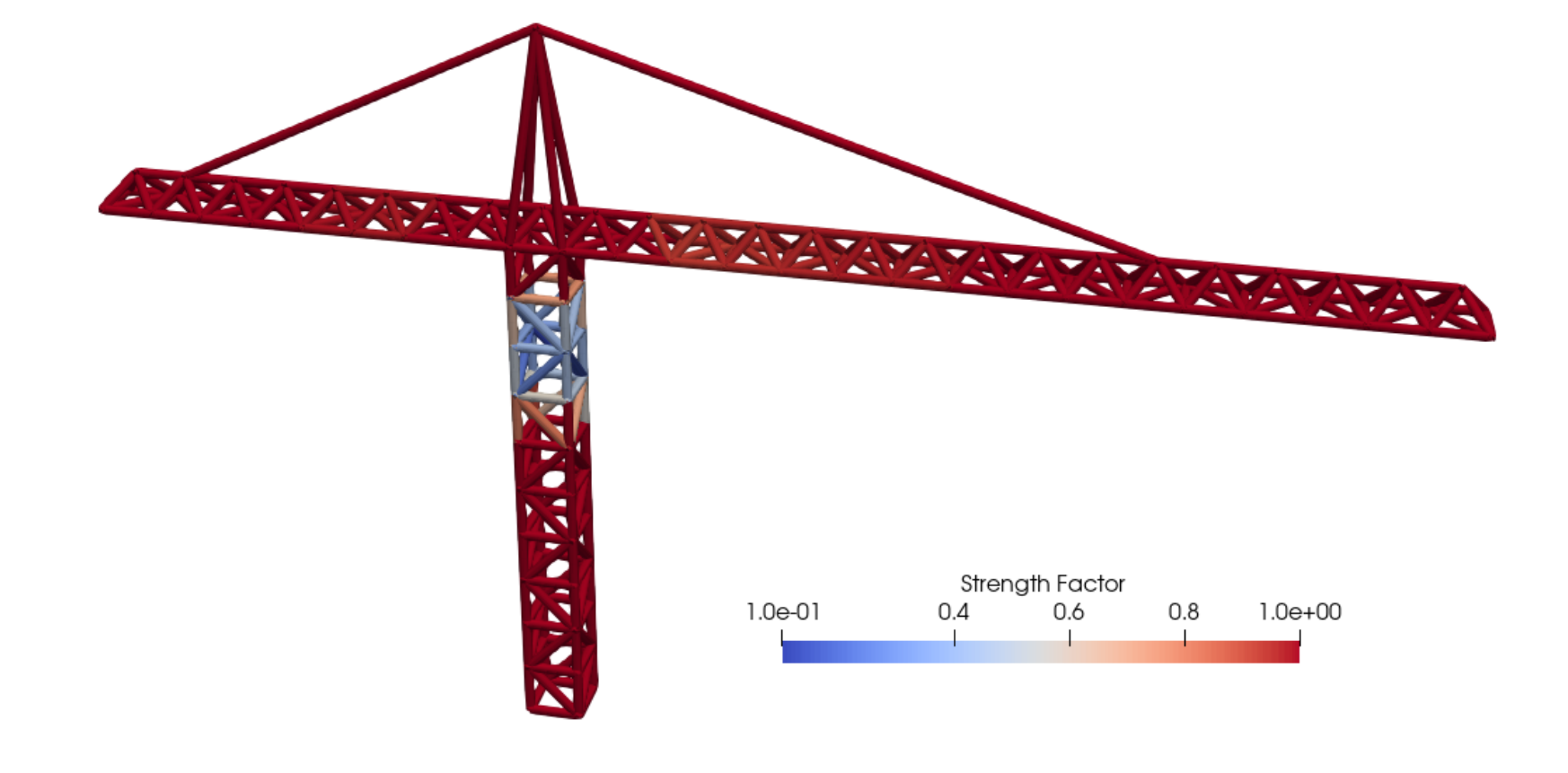}} 
    \caption{Top row:
(Left) target strength factor and sensor location. 
(Right) Weakness identification with standard expectation in the objective function. 
Bottom row: Identification using a risk-measure such as CVaR$_\beta$ with different confidence levels $\beta = 0.3$ (left)
and $\beta = 0.8$ (right).}
    \label{f:CVaR}
\end{figure}

In case the load $\bm f$ is uncertain, the results are shown in Figure~\ref{f:CVaR} which compares 
the results obtained from using standard expectation and a risk-measure, i.e., conditional value at 
risk (CVaR) with a confidence level given by $\beta$ \cite{FAiraudo_HAntil_RLoehner_URakhimov_2024a}. 

\subsection{Medical Digital Twins}

President Obama announced his precision medicine initiative during the 2015 state of the union address and a working group was established \cite{hudson2015precision}. 
The plan was to collect data on a million patients to create a framework leading to individual patient specific medical treatment. Digital twins from engineering have a strong analog in medicine as a tool for truly personalizing medical care. First attempts at medical digital twins (MDTs) date back two decades \cite{eddy2003archimedes}. 
Since then the literature has grown significantly. The current state of the art is summarized in \cite{RLaubenbacher_2024a}.

There are relatively few examples of models in medicine that directly fit the digital twin paradigm, as laid out during the workshop. An example that comes closest is the artificial pancreas that provides almost automated insulin delivery for type I diabetes patients. There are several examples in cardiology, as well as in cancer. Specific examples and references are included in \cite{RLaubenbacher_2024a}.

\newpage
\section{Summary of Funding Agency Panel}

This panel was led by Dr. Yuliya Gorb from NSF. The panelists were:
Jodi Mead (NSF), Fariba Fahroo (AFOSR), Warren Adams (AFOSR),
Tim Bentley (ONR), Meredith Hutchinson (ONR), Robert Martin (ARO),
Yannis Kevrekidis (DARPA), Ali Ghassemian (DOE), Grace Peng (NIH).

\noindent
The following questions were prepared and were sent to all the panelists 
in advance. 
\begin{enumerate}
\item What specific initiatives or funding opportunities do all participating agencies offer to support research activities focused on advancing the mathematical foundations of Digital Twins technology, if any? Overall, can your funding agency prioritize research projects focused on developing novel mathematical models and/or algorithms to enhance the development of Digital Twins technology? If not, what specific areas of research does your agency identify as necessary to facilitate the development or advancement of Digital Twins technology?
 
\item What level of emphasis does your funding agency place on supporting research exploring the theoretical aspects of Digital Twins in mathematics, if any, compared to those focusing on practical implementation? Could you specify the particular mathematical challenges your agency expects researchers to explore in this field? Could you outline the criteria or key areas of interest that your funding agency will consider when evaluating proposals related to Digital Twins in mathematics, if any?
 
\item Could you provide examples of successful research projects potentially funded by your agency that could significantly contribute to the advancement of mathematics within the realm of Digital Twins?
 
\item Are there any collaborative funding opportunities or partnerships encouraged by the government overall to promote interdisciplinary research integrating mathematics with Digital Twins applications in various fields?
 
\item How does your funding agency ensure the dissemination of findings/knowledge obtained/acquired through funded research projects in domains of ‘Digital Twins in Mathematics’ or ‘Mathematics of Digital Twins’, if any, to benefit academia, industry, and society at large?
\end{enumerate}

\noindent
{\bf Response summary:}
Overall there was a significant enthusiasm from various funding agencies to support 
MATH-DT. There were several logistical challenges that were recognized, which could 
potentially slow the DT research progress.
\begin{enumerate}
	\item There are some recent examples of individual efforts from the funding agencies to support 
		MATH-DT related research. For instance, a recent MURI from AFOSR and a project from ARO 
		in Europe. Several examples of DTs were discussed, in particular, power system under DOE. 

	\item It was discussed that different project types may be funded in the future calls:
		\begin{itemize}
			\item Smaller projects consisting of at least one domain expert and one mathematician.
			
			\item Larger projects with larger interdisciplinary teams, with mathematics at their core. 
				
			\item MATH-DT research institutes, something similar to recent NSF funded institutes, 
				could be highly beneficial. Examples of NSF funded institutes are: National 
				AI Research Institutes, ICERM (Brown University), Institute for Mathematical and 
				Statistical Innovation (Chicago).
		\end{itemize}
	
	\item It was also discussed that several of the agencies typically tend to support 
		more basic research and have smaller research projects. However DT  
		research require researchers working in both basic sciences and domain 
		specific science. Can agencies come together to support this?
		
	\item A need to engage in a partnership between academia, national labs and industry was
		also discussed.
	
	\item All the agencies acknowledged that innovative ideas are needed to train next generation
		scientists who are interested in DTs.
\end{enumerate}

\section{Dream Big Session}
\label{s:dream}

DB1-DB5 were the five topics that were identified for the Dream Big sessions. The workshop 
attendees divided themselves in various groups based on the topic of their choice. In total, 
there were 7 groups: DB1 (1), DB2 (2), DB3 (2), DB4 (1), DB5 (1), see Table~\ref{t:DB}. 
\begin{table}[htb!]
\centering
\begin{tabular}{l|l} \hline 
{\bf DB1: Data (1)}  		    &  {\bf DB2: Modeling \& Forward Problem (2)}   \\ \hline
Acquisition of dynamic data   		    & Data assimilation				 	   \\
Data curation and management 	    & Multiphysics and multiscale  		  	    \\
Big data						    & Uncertainty						    \\ 	
Data assimilation		   	    	    & Surrogates and reduced order models  	    \\ 
DNNs and ML			  		    & Sensor placement 		 			    \\
Sensors					    	    & Optimal experimental design		            \\
				   		    	    & Scientific ML						    \\ 
\hline \hline
{\bf DB3: Optimization \& Inverse Problems (2) }  &  {\bf DB4: Validate, Verify \& Predict (1)}  \\ \hline
Optimal control \& optimal design	    	    &  Validation of DT 			\\ 
Nonconvex, nonsmooth, 	\& multi-objective  &  Update of DT 				\\
Uncertainty 					   	    &  Verification of DT   			\\
Surrogates and reduced order models	    &	UQ						\\
Data assimilation					    & 	Maintenance				\\	
Sensor placement  					    & 	Error \& sensitivity analysis 	    	      	       			\\ 
Optimal experimental design			    &						        \\  
Scientific ML 						    & 							\\ \hline\hline
\end{tabular}	
\begin{tabular}{||c||}
{\bf DB5: Software (1)} \\ \hline\hline
\end{tabular}
\caption{\label{t:DB}Dream Big (DB) session was divided into five groups listed above. 
Each group had further sub-groups as indicated by the numbers in the parenthesis.}
\end{table}
On the first day, each group discussed DB1-DB5 in the context of Q1-Q3. On the second 
day, each group presented their findings to all the conference participants to receive a 
feedback from all the attendees. On the final day, each group again met and discussed 
the previous comments and, refined their key points. 
Detailed comments from each group can be found in Appendix~\ref{app:DB}. The key findings are summarized in the next section.

\newpage

\section{Overall Summary}
\label{s:recom}

\subsection{Summary of Challenges}
\label{s:chal}

\begin{enumerate}\itemsep0pt

\item {\bf Better Pre-Processors}: As explained before, DTs require 
constant updating of the physical parameters, boundary conditions,
loads and in some cases even geometry for the physical models that
describe them. This implies the need of updating or `curating' the
DT models - for each of the disciplines (structural-, fluid-,
electro-, thermo-, etc. dynamics). And once a change occurs (e.g. a
part is changed in a turbine or airplane), the complete set of models
need to be re-run to make sure the product/process is safe and the
DT is updated. This process, if not automated, will be extremely expensive,
hampering the widespread adoption of DTs. Given that most products and 
processes nowadays are first computed, then built/executed, it is therefore
imperative to develop computer aided design (CAD) tools and m\'etier-specific
pre-processors that allow an automatic update of DTs. These would then
yield `DT-Ready' models that can be updated automatically throughout the
lifecycle of the product/process. 

In a broad sense, we need DTs which can incorporate data and we need to 
co-design softwares and algorithms which are tuned to DT framework.

\item {\bf Gaps In Forward Problems}: In many areas of engineering, biology,
medicine and material science the basic physics of the problem are still
unknown. This implies that at this point we do not know which are the
relevant PDEs or ODEs that properly describe the physics. In other instances,
we know the basic physics, but the uncertainty in the physical parameters is
very large (this could be due to material parameter variations or the 
inability to conduct experiments/obtain good measurements). We are aware that
the research required here does not fall into the realm of mathematics per se. 
However, as new PDEs or ODEs that describe the newly found physics emerge, 
it is important for math to obtain rigorous understanding of the basic properties 
of these PDEs and ODEs, such a uniqueness, stability, convergence, asymptotic 
behaviors (e.g. formation of discontinuities), and reduced order modeling. 

\item {\bf Coupling of Models, Physics, Scales, ...}: For whole classes of
problems, the description of the physical reality even with simple models 
may require a coupling of models (fluid, structure, heat, ...) across many
length (atomistic/continuum) and time (picoseconds to hours) scales.
Other problems may require the coupling of models of different dimensional
abstraction (e.g. beams, shells and solid elements in a finite element
model of a bridge). Mathematics plays a key role in assuring that this
coupling of models, physics and scales leads to unique, stable and
convergent answers. Otherwise, no valid DT can be built.

\item {\bf Optimization and Inverse Problems}: In order to infer properties of the real-life
object/process of the DT, inverse problems need to be solved. 
Similar scenario occurs when a DT assists a human in 
making a decision about the physical model. The most expedient 
way to obtain for instance the unknowns is using optimization via
the so-called method of adjoints. The theory of 
adjoint solvers is well developed for classic PDEs like elasticity or
incompressible flow, but for more complex physics/models, or models
that are the result of a coupling of models, the adjoints may be
non-existent, difficult or impossible to obtain. Even in the classical 
settings, such as bilevel optimization, optimal experimental design 
or situations with discrete design variables, the situation is extremely  
challenging. These problems are nonconvex, nonsmooth and 
may also be multi-objective.
This implies that there is a real need for further in-depth 
mathematical analysis and algorithm development for these 
`difficult problems and adjoints'. 

\item {\bf Uncertainty Quantification and Optimization Under Uncertainty}:
Uncertainty should be taken into account in every step of the decision process
whether it is during modeling, forward problems, optimization or engineering
design. Only then can one make meaningful decisions. For example, one 
must account for uncertainty and extreme case scenarios, for instance 
using, risk-measures. The resulting forward and optimization problems 
with uncertainty are nonconvex, nonsmooth and high dimensional and 
cannot be solved using traditional approaches. One should also take 
into account that the underlying probability distribution maybe unknown. 

\item {\bf Software Development}: While designing softwares to carry out 
efficient adjoint calculations it is important to take into account both
the direct methods and automatic differentiation, by developing integrated 
approaches allowing modularity and co-designing methods and software. 
See the first bullet above and DB5 in Appendix~\ref{app:DB}. 

\item {\bf Workforce Development}: DTs, by their very nature, are the
result of multi-disciplinary effort. They involve model preparation,
model execution, sensor input, system output analysis, inferring
the state of the system from measurements, and model update/upkeep.
In order to effectively deploy DTs and derive their potential benefits
for society (increased safety, increased comfort, longer life cycles,
reduced environmental footprint) the current and future developers 
and users of DTs need to be trained properly. This mirrors similar
requirements that first arose when computational mechanics began to
have an impact in product/process design and analysis. For DTs, the
workforce will need a combination of classic and modern numerical methods 
(PDEs, optimization, numerical linear algebra, randomized methods, 
reduced order modeling), 
statistics and uncertainty quantification (sensors, forces, actuators), 
computer science (basic programming skills, large data sets), and, of
course, engineering.

\end{enumerate}

\subsection{Summary of Opportunities}
\label{s:opp}

\begin{enumerate}				
	\item {\bf Canonical Testcases}: In every field of science, numerical methods
are gauged and verified against classic or `canonical' testcases. In 
aerodynamics it is the NACA0012 airfoil, in hydrodynamics the Wigley Hull,
in structural mechanics the plate with a hole, in ML it is MNIST/CIFAR dataset, etc. 
In order to develop 
robust DTs, a series of canonical testcases that are open, accessible and
verifiable need to be defined. Ideally, they should be backed up by real-life
instantiations. In all of these cases the numerical model(s) of the real-life
object would be accessible in a repository, with several copies of varying
detail and size possible. 

{\bf For Civil Engineering:} this could be a simple bridge, where measurements are
taken at several locations, and certain trusses, beams or plates are
artificially weakened in order to see if the DT correctly predicts the
location.

For Mechanical Engineering: this could be a simple vibrating beam with an
attached mass. As before, measurements (either displacements, strains or
vibration data (displacement/frequencies)) are taken at several locations, 
and some part is weakened (or removed) in order to see if the DT correctly 
predicts the location.

{\bf For Biomedical Problems:} The most common strategy for solving biomedical problems, such as the discovery of mechanisms underlying biological dysregulation and therapeutic approaches to correct or prevent it, is to use \emph{in vitro}, \emph{in vivo}, and \emph{ex vivo} experimental models. \emph{In vitro} models include traditional 2D cell cultures of human or animal cell lines or primary cells derived from donors. In recent years, 3D cell cultures, organoids, organ-on-a-chip systems, and other platforms have been developed. Tissue culture systems that allow preservation and manipulation of 3D intact human and other tissue have been developed that allow time course experiments. Finally, \emph{in vivo} animal models include mouse, rat, pigs, nonhuman primates, and others. 

One strategy to develop a suite of canonical test cases is to choose several biomedical use cases and develop experimental models for them at different levels of complexity. The experimental models should be such that one can set up the dynamic interplay between digital and physical twin. This requires taking repeated measurements on the physical twin to update the model and perform interventions on the physical twin obtained from interrogating the digital twin. 

An expanded version of this strategy could begin with a computational model that plays the role of the physical twin and can be ``learned" from repeated observations from the model. 

{\bf For Medicine:} In order to develop canonical test cases that relate directly to patients, it is important to choose a collection of clinical conditions that enable relatively easy collection of large and heterogeneous time course data sets from patients. An example is sepsis, treated in the ICU, which provides the opportunity to develop such data sets. For instance, one could leverage existing efforts, such as in \cite{sepsis}.

\item {\bf Interdisciplinary Cross-Pollination}:
As stated before, DTs by their very nature require a multidisciplinary
approach in order to make progress. This opens the possibility of
interdisciplinary cross-pollination, where successful solution approaches
are migrated between fields. These could be (and this is a non-exhaustive
list) in the areas of inverse problems, optimization, field solvers 
(forward and adjoint), multiscale solution approaches, and multifield 
solution approaches.

 The multidisciplinary approach that is suggested represents a unique opportunity to develop interdisciplinary training and education.

\item {\bf Hierarchy of Digital Twins}:
DTs can mean different things to different users or different fields.
Let us take the case of commercial airplanes. The airplane manufacturer
can follow each of the planes as a distinct DT. The models that were used
during preliminary and detailed design and analysis are available, and
preceded manufacturing. A given airplane leaves the factory and starts 
flying. Maintenance is required every so often, and parts may be exhanged
(the usual wear and tear). At the most abstract level, the DT is simply
`plane x', and the data the `sensors' send is simply: hours flown,
passenger/cargo carried, parts changed (in this case the logbook is the
`sensors'). At the next level, one could envision the enumeration of all
the parts that comprise the airplane: beams, trusses, motors, seats, ...
The sensors for such an `inventory only' DT would again be logbooks,
purchase orders, maintenance certifications, etc. One could add more
and more details, until one reaches a finite element discretization
of the complete plane. The sensor data is now displacement or strain 
data at several locations on the airplane. And from this data one
infers the state/health/safety of the structure. There could be several
levels of structural models: some would have a higher level of abstraction 
(e.g. beams, lumped masses) than the more detailed ones (e.g. with more
plates, shells and solid elements).
Note that a hierarchy of models such as the one outlined above is 
encountered in any field, be it biology, engineering, medicine, climate 
modeling, etc.

\item {\bf Characterization of Rare Events}:
Accidents happen when unforeseen or rare events take place. Given that
DTs have to work with uncertain data (models, sensors, loads), together 
with risk-measures, they are uniquely capable of characterizing rare events. 
This should have a considerable impact on safety for products and processes 
that have DTs.	
\end{enumerate}

\newpage
\bibliographystyle{plain}
\bibliography{refs}

\newpage

\appendix

\section{Dream Big}
\label{app:DB}

This section provides detailed notes from the discussion during the Dream Big sessions.
An overall summary is provided in section~\ref{s:recom}. 

\subsection{DB1: Data}

This topic is covered throughout DB2-DB5.

\subsection{DB2: Modeling \& Forward Problem}

Modeling and forward problems form the foundation of DTs. Because of their multi-scale, 
multi-physics nature, and additional uncertainty in model (or data), many challenges and 
opportunities arise in DTs. In particular:

\bigskip 
\noindent
{\bf Group 1:} Led by Akil Narayan (The University of Utah) 
\paragraph{Challenges and Opportunities}
\begin{enumerate}\itemsep0pt
	\item {\bf Coupling}
		\begin{itemize}
			\item \emph{Challenges:} To tackle DT problems one has to face the reality of 
			heterogeneous models with multi-scale, multi-physics features and multiple
			components. 			
				
			\item \emph{Opportunities:} Develop impactful, explicit paradigms. In particular, it is 
				critical to study stability and robustness at coupling interfaces. 
				It is also interesting to create robust bridges between actual DT instantiation and 
				mathematical model.
		\end{itemize}
	
	\item {\bf Trustworthiness of ML models}
		\begin{itemize}
			\item \emph{Challenges:} Explainability, interpretability, certifiability are challenging for DTs. 
				ML can be flexible, but is it trustworthy, especially compared to ``traditional" numerical methods resulting from 
				decades-long development?
				
			\item \emph{Opportunities:} Develop flexible, efficient and general purpose forward solvers, 
				by using inexactness, adaptivity and leveraging interdisciplinary connections. 
		\end{itemize}
		
	\item {\bf Model uncertainty}
		\begin{itemize}
			\item \emph{Challenges:} Where does model uncertainty belong? With mathematicians, 
				modelers, engineers, etc?
				
			\item \emph{Opportunities:} Develop frameworks for codifying model uncertainty, and create 
				more precise understanding of ``unknown unknowns", for instance, unknown probability
				distribution. 
		\end{itemize}	
		
	\item {\bf Workforce development}
		\begin{itemize}
			\item \emph{Challenges:} Hard to train early career researchers to focus on DT’s. 
				Incorporate formal training, hands-on experience.
				
			\item \emph{Opportunities:} There are possibilities to create training paradigms. 
				Funding agencies need to play a leading role. 
		\end{itemize}			
\end{enumerate}

\bigskip 
\noindent
{\bf Group 2:} Led by Rainald L\"ohner (George Mason University). Several of the suggestions and 
recommendations from this group have been directly incorporated in section~\ref{s:recom}. 

\subsection{DB3: Optimization \& Inverse Problems}

Optimization and inverse problems arise in many facets of DTs. For example, the 
integration of data into the digital twin can be formulated as data assimilation or 
inverse problem. Referring back to Section~\ref{s:examples}, several examples of digital 
twins include a feedback from twin to data that involves optimal control or the optimal 
design of experiments to optimally  update the digital twin. The resulting optimization problems 
are high-dimensional, discrete, multi-modal, involve uncertainties, \dots The following 
challenges and opportunities were identified by the two groups:

\bigskip 
\noindent
{\bf Group 1:} Led by Sven Leyffer (Argonne National Laboratory) 
\paragraph{Challenges}
\begin{enumerate}\itemsep0pt
	\item {\bf Risk quantification and risk mitigation for decisions:} For example, rare events and optimization with chance constraints or optimization under uncertainty. Results in nonconvex large-scale optimization problems with high-dimensional uncertainty. How do we provide performance guarantees for DTs? 
	
	\item {\bf Design of experiments with non-Gaussian errors:} Additional challenges include correlated observations, model errors, and nonlinear forward models, integrated within a larger feedback loop (DT).
	
	\item {\bf Integrating initial design, control, and experimental design:} Life-cycle optimization; different optimization problems. Results in bilevel or hierarchical optimization problems.
	
	\item {\bf Bilevel or Hierarchical Optimization:} Bilevel, hierarchical optimization with nonconvex lower-level problems (for example, Stackelberg games). No clear reformulation as a single-level optimization (max=min otherwise), hence not clear how to get a tractable formulation.
	
	\item {\bf Discrete design or control variables:} E.g. branching logic (controls or states) results in variational inequalities (e.g. ReLU). Problems have distributed controls or designs in infinite dimensions. Lack of optimality conditions.
\end{enumerate}

\paragraph{Opportunities}\itemsep0pt
\begin{enumerate}
	\item Ability to exploit hierarchy of models in DTs for design, control, experiments.
	\item Life-cycle analysis: cheaper, more efficient, more reliable and robust.
	\item Enable new discoveries, by running experiments in-silico.
	\item Can we develop rigorous AI for science with error analysis by building on DTs?		
\end{enumerate}

\bigskip 
\noindent
{\bf Group 2:} Led by Bart van Bloemen Waanders (Sandia National Laboratory) 
\paragraph{Challenges and Opportunities}
\begin{enumerate}\itemsep0pt
	\item {\bf Multi-physics/scale/fidelity/components} 
	\begin{itemize}\itemsep0pt
		\item adjoint, Hessians, higher order derivatives, priors
		\item interface of components
		\item modeling error, parametric uncertainties
		\item high dimensionality of state, optimization, uncertainty: computational efficiency
		\item Incorporation of reduced order models (ROMs)
	\end{itemize}
	\item {\bf Beyond academic model problems}
	
		\begin{itemize}\itemsep0pt
			\item benchmark digital twin model datasets – bridge, fusion experimental facilities 
		\end{itemize}		
		
	\item {\bf Toolbox accessibility} 
		\begin{itemize}
			\item documentation
			
			\item fuse modeling, optimization, UQ, decision-making
		\end{itemize}
		
	\item {\bf Mathematical advances}		
		\begin{itemize}
			\item Mixed integer programming 
			\item Optimization under uncertainty
			\item goal-oriented data acquisition
		\end{itemize}
		
	\item {\bf Training of next generation researchers}
	
		\begin{itemize}\itemsep0pt
			\item inter-disciplinary programs in academia, different incentive structure
			\item industrial, laboratory, academic partnerships
			\item Integrate domain expertise
		\end{itemize}		
\end{enumerate}

\subsection{DB4: Validate, Verify \& Predict}

\bigskip 
\noindent
Led by Ansu Chatterjee  (University of Maryland Baltimore County)

\paragraph{Opportunities and Challenges}
\begin{enumerate}\itemsep0pt
	\item \emph{Opportunity:}			
		{\bf Balance short-term and long-term uncertainties}  in the DT according to its use 
		and utility. Control uncertainty growth rates in specific dimensions. 
		
		\emph{Challenge:} which uncertainties to focus on and control and how?
	
	\item {\bf Process control, robustness and quality maintenance of DT:} 
	
		\emph{Opportunity:} 
		decide on the frequency of collecting data, frequency of updates, and whether the DT 
		as an industrial/clinical process is ``in control".
		
		\emph{Challenge:} how to measure and communicate the 
		quality of a DT? classical statistical and machine learning assumptions do not hold. What are 
		the right quantities to measure? 

	\item {\bf Design DTs with upfront uncertainty controls and for multiple competing goals:} 
	
		\emph{Opportunity:} 
		(like safety and efficiency) by allocating limited resources optimally: multi-armed bandit 
		approaches, exploration-exploitation.
		
		\emph{Challenge:} adaptive 
		controls is hard, classical assumptions do not hold, (reinforcement) learning with complex 
		structure/data is open.

	\item {\bf Track internal variability of the physical system:} the variability in a built environment is not 
		the same as that of a society. 

	\item {\bf Account for data quality issues and ethics:} 
	
		\emph{Opportunity:} data from the physical system may be restricted 
		due to cost, ethics, technology, may have missing values, may have selection and other biases. 
		
		\emph{Challenge:} biased data, (informative) missing values 
		and sampling issues. 
		
	\item {\bf Visualization of uncertainties:} 
	
		\emph{Opportunity:} efficient and honest communication of uncertainties. Design new tools
		to communite uncertainties.
	
		\emph{Challenge:} Gaps in education and general understanding of uncertainties.	
\end{enumerate}

\subsection{DB5: Software}
\label{s:soft}

\bigskip 
\noindent
Led by Carol S. Woodward  (Lawerence Livermore National Lab)

\paragraph{Challenges and Opportunities}
\begin{enumerate}
	\item Efficient adjoint calculations through the numerical and application software stack
		\begin{itemize}
			\item Automatic differentiation (AD)
		\end{itemize}
		
	\item Integrated methods that allow for modularity in software implementations
		\begin{itemize}
			\item Need to develop methods that are aware of and can exploit the reason for their need 
				and other methods using them and used by them
			\item E.g., integrate uncertainties across the full numerical stack while allowing for flexibility 
				in specific methods for computing uncertainties
		\end{itemize}
		
	\item Co-Design of methods and software 	
		\begin{itemize}
			\item Data acquisition mechanisms, software that includes new data, and mathematical 
				methods that can adapt the system structure based on new data
				
			\item Research how to integrate workflows into the full methods stack (forward simulation, 
				optimization, solvers)
				
			\item Visualization of data and results along with calculations of quantities of interest, including 
				uncertainties					
		\end{itemize}
		
	\item Software hardening and support; incentivize new projects to reuse already made products
		\begin{itemize}
			\item Software hardening - documentation, testing, community building (tutorials, user 
				communities), …
			
			\item Method documentation beyond core parts using Appendices or Supplemental Materials 
				sections of publications (do we want these refereed?)
				
			\item Need to have a place to host reusable software
			
			\item Require funding for Research Software Engineering (RSE) for any project whose goal 
				is to develop a digital twin
				
			\item Need to find ways to value the work of software hardening and RSE
				\begin{itemize}
					\item Can we value User Guides or Tutorials as publications?
					
					\item What is the career path for an RSE
				\end{itemize}
				
			\item Open source development is viewed as very critical for helping with this (community 
				building, community contributions, etc.)
		\end{itemize}
		
	\item Curriculum development to include numerical software development and good software design 
		and implementation practices, some RSE skills
		\begin{itemize}
			\item Debugging skills
			\item Performance analysis
		\end{itemize}
\end{enumerate}

\newpage

\section{Schedule}
\label{s:schedule}

\begin{center}
{\large\bf Schedule of Mathematical Opportunities in Digital Twins (MATH-DT)} \\
December 11--13, 2023\\
\end{center}

\underline{Monday, December 11, 2023}

\medskip
\begin{tabular}{ll}
8:00AM - 9:00AM  & Breakfast / Registration \\
8:45AM - 9:30AM  & Welcome address by the Organizers, Andre Marshall (VP Research, \\ 
			     &  GMU), and David Manderscheid (NSF DMS Division Director)  \\
9:30AM - 11:00AM & \cellcolor{red!25}{Session 1 (Plenary Talks)} \\
11:00AM - 11:30AM  & Coffee \\
11:30AM - 12:30PM  & \cellcolor{blue!25}{Session 2 (Applications Panel)} \\
12:30PM - 1:30PM  & Boxed Lunch \\
1:30PM - 3:00PM   & \cellcolor{red!25}{Session 3 (Plenary Talks)} \\
3:00PM - 4:00PM   & \cellcolor{cyan!25}{Dream Big \& Report - Discussion} \\
4:00PM - 5:30PM   & \cellcolor{yellow!80}{Poster / Demo Session} 
                  (Jointly with SIAM DC-Baltimore Section)\\
5:30PM - 6:00PM   & \cellcolor{yellow!80}{SIEMENS (Demo)}           
\end{tabular}

\medskip
\underline{Tuesday, December 12, 2023}

\medskip
\begin{tabular}{ll}
8:00AM - 9:00AM  & Breakfast / Registration \\
9:00AM - 10:30AM & \cellcolor{red!25}{Session 4 (Plenary Talks)} \\
10:30AM - 11:00AM  & Coffee \\
11:00AM - 12:00PM  & \cellcolor{blue!25}{Session 5 (Applications Panel)} \\
12:00PM - 1:00PM  & Boxed Lunch \\
1:00PM - 2:30PM   & \cellcolor{green!25}{Session 6 (Funding Agency Panel)} \\
          		        & (Jointly with SIAM DC-Baltimore Section)\\			     	
2:30PM - 4:30PM   & \cellcolor{cyan!25}{Dream Big \& Report II}  \\                
4:30PM - 5:00PM  & \cellcolor{yellow!80}{NVIDIA (Demo)} \\                  
5:30PM -  7:00PM        & Conference Dinner 
\end{tabular}

\medskip
\underline{Wednesday, December 13, 2023}

\medskip
\begin{tabular}{ll}
8:00AM - 9:00AM  & Breakfast / Registration \\
9:00AM - 10:30AM & \cellcolor{red!25}{Session 7 (Plenary Talks)} \\
10:30AM - 11:00AM  & Coffee \\
11:00AM - 12:15PM  & \cellcolor{cyan!25}{Dream Big \& Report III - Final Presentations} \\
12:30PM - 1:30PM    & Boxed Lunch
\end{tabular}

\noindent
{\bf Notes}
\begin{enumerate}\itemsep0pt
    \item MATH-DT has ``in Cooperation with Association for Women in Mathematics (AWM)" status. This conference supports the Welcoming Environment Statement of the AWM.

    \item Each {\bf Plenary talk} will be 20-minutes followed by a moderated panel discussion of 30 minutes. 

    \item Each {\bf Applications Panelist} will give 10-minute talk followed by a moderated panel discussion of 30-minutes.

    \item Each {\bf Funding Agency} will give 7-8-minute introduction followed by a moderated panel discussion of 30-minute.

    \item Session assignments are on the next page.
\end{enumerate}


\section{Speaker List}
\label{s:speakers}

\noindent
{\bf Session 1: Plenary Talks}\\
Omar Ghattas
(The University of Texas Austin, USA) \\
Julianne Chung
(Emory University, USA)\\
Reinhard Laubenbacher
(University of Florida, USA)

\bigskip
\bigskip
\noindent
{\bf Session 2: Applications Panel I}\\
Stefan Boschert (SIEMENS) \\
Sharon Di (Columbia University)\\
Mustafa Megahed (ESI Group)

\bigskip
\bigskip
\noindent
{\bf Session 3: Plenary Talks}\\
Roland W\"uchner
(Technical University Braunschweig, Germany)\\
Rainald L\"ohner (George Mason University)\\
Bart Van Bloemen Waanders
(Sandia National Laboratories) 

\bigskip
\bigskip
\noindent
{\bf Session 4: Plenary Talks}\\
Enrique Zuazua
(Friedrich-Alexander-Universit\"at Erlangen-Nürnberg, Germany)\\
Carol Woodward
(Lawerence Livermore National Laboratory, USA)\\
Karen Veroy-Grepl
(Eindhoven University of Technology, Netherlands)

\bigskip
\bigskip
\noindent
{\bf Session 5: Applications Panel II} \\
Brent Bartlett (NVIDIA)\\
Aldo G Badano (FDA)\\
Juan Cebral (George Mason University)

\bigskip
\bigskip
\noindent
{\bf Session 6: Funding Agencies Panel} 

\medskip
\begin{tabular}{|l|l|}\hline 
National Science Foundation & Yulia Gorb and Jodi Mead \\
Air Force Office of Scientific Research & Fariba Fahroo and Warren Adams \\
Office of Naval Research & Tim Bentley and Meredith Hutchinson \\
Army Research Office & Robert Martin  \\
DARPA  &  Yannis Kevrekidis \\
Department of Energy &  Alireza Ghassemian \\
National Institute of Health & Grace Peng \\ \hline
\end{tabular} 

\bigskip
\bigskip
\noindent
{\bf Session 7: Plenary Talks}\\
Arvind K. Saibaba
(North Carolina State University, USA) \\
Akil Narayan
(University of Utah, USA)\\
Lise Marie Imbert-G\'{e}rard
(University of Arizona)

\newpage

\section{Titles and Abstracts}
\label{s:abstracts}

\begin{center}
 \underline{\bf\Large Plenary Speakers}
\end{center}

\vspace{1cm}

\begin{enumerate}[1.]
\setlength\itemsep{2em}

\item 
\HA{\large Julianne Chung} (Emory University)

\smallskip
\noindent
{\bf Title.} Feedback Flow from Physical to Virtual: Mathematical Opportunities in Inverse Problems, Uncertainty Quantification, and Data Assimilation

\smallskip
\noindent
{\bf Abstract.} 
Digital twins have gained widespread interest recently as a tool for accelerating scientific discovery and revolutionizing industries. Central to any successful digital twin is the bidirectional interaction between the virtual representation and the physical counterpart.  In this talk, we focus on the feedback flow from physical to virtual, that is, the process of combining physical observations to update the virtual models in a rigorous, systematic, and scalable way. We describe challenges and opportunities for digital twins in the field of inverse problem methodologies, uncertainty quantification, and data assimilation.

\item 
\HA{\large Omar Ghattas} (The University of Texas)

\smallskip
\noindent
{\bf Title.} Geometric Neural Operators for Digital Twins 

\smallskip
\noindent
{\bf Abstract.} 
A digital twin is a computational model of a physical system that
continually updates its knowledge of the system by assimilating
observational data, and in turn informs decisions and controls the
system to achieve a desired goal, over a continually evolving time
horizon. When the model takes the form of a system of PDEs, and when
the data assimilation and optimal control problems must be solved
under uncertainty, the computations underlying digital twins quickly
become prohibitive.  Efficient surrogates of the maps from parameters
or decision variables to objective functions, which involve solution
of the forward PDEs, is key to making digital twins tractable.
Unfortunately, constructing such surrogates presents significant
challenges when the parameter/decision dimension is high and the
forward model is expensive, which limits the amount of available data.

Deep neural networks have emerged as leading contenders for overcoming
the challenges of constructing infinite dimensional surrogates, which
are known as neural operators. We demonstrate that black box
application of DNNs for problems with infinite dimensional parameter
fields leads to poor results when training data are limited due to the
expense of the model. Instead, by constructing a network architecture
that captures the geometry of the map -- in particular its smoothness,
anisotropy, and intrinsic low-dimensionality -- one can construct a
dimension-independent reduced basis neural operator with superior
generalization properties using just limited training data. We employ
this reduced basis neural operator to make tractable the solution of
PDE-constrained Bayesian inverse problems and optimal control under
uncertainty, with applications to inverse wave scattering, inverse
hyperelasticity, inverse earthquake subduction, and optimal flow
control.

This work is joint with Thorsten Becker, Lianghao Cao, Peng Chen,
Dingcheng Luo, J. Tinsley Oden, Thomas O'Leary-Roseberry, Simone Puel,
Umberto Villa, and Keyi Wu.

\item 
\HA{\large Lise Marie Imbert-G\'{e}rard} (University of Arizona)

\smallskip
\noindent
{\bf Title.} Nuclear Fusion reactors \& Digital Twins

\smallskip
\noindent
{\bf Abstract.} 
Nuclear fusion could become a clean, safe and abundant source of energy. Several techniques are being explored by the scientific community to confine matter in the plasma state and control nuclear fusion, including magnetic confinement fusion, leveraging an external magnetic field to confine charged particles. Modeling for magnetic confinement fusion raises many challenges: the very nature of controlled fusion is inherently multi-scale; the coupling between plasma and electromagnetic fields is inherently multi-physics; the complexity of efficient plasma modeling is inherently multi-fidelity. Moreover, vast amounts of data are generated and collected in experimental fusion reactors, throughout their lifespan, from design to construction and operation. In this context, developing digital twins of reactors can provide a unique perspective on the pursuit of nuclear fusion.
This talk will present a brief review of existing efforts in this direction, and discuss related challenges and opportunities.

\item 
\HA{\large George Karniadakis} (Brown University)

\smallskip
\noindent
{\bf Title.} Modular Building of Digital Twins using Neural Operators

\smallskip
\noindent
{\bf Abstract.} 
We will present a new data assimilation framework, the DeepM\&Mnet, for simulating multiphysics and multiscale problems at speeds much faster than standard numerical methods using pre-trained neural networks (NNs). We first pre-train DeepONets that can predict independently each field, given general inputs from the rest of the fields of the coupled system. DeepONets can approximate nonlinear operators and are composed of two sub-networks, a branch net for the input fields and a trunk net for the locations of the output field. DeepONets, which are extremely fast, are used as building blocks in the DeepM\&Mnet and form constraints for the multiphysics solution along with some sparse available measurements of any of the fields. We demonstrate the new methodology and document the accuracy of each individual DeepONet, and subsequently we present two different DeepM\&Mnet architectures that infer accurately and efficiently fields in electroconvection, hypersonics and oil \& gas applications.  The DeepM\&Mnet framework is general and can be applied for building any complex multiphysics and multiscale models based on very few measurements using pre-trained DeepONets in a plug-and-play mode.

\item 
\HA{\large Reinhard Laubenbacher} (University of Florida)

\smallskip
\noindent
{\bf Title.} The Mathematics of Medical Digital Twins

\smallskip
\noindent
{\bf Abstract.} 
The digital twin concept has its origins in industry. One industrial
equipment manufacturer advertised its digital twin capabilities to its customers as 
"No unplanned downtime" for its products. There is a compelling aspirational analog in 
healthcare: 'No unplanned doctor visits." Of course, the challenges of building
digital twins for human patients are incomparably greater
than for machinery. Nonetheless, there are now several instances of what 
might be called digital twins in medicine, and many more ongoing development
projects. Aside from our incomplete understanding of human biology, 
relative sparseness of data characterizing human patients, and logistical 
difficulties in implementing computational models in healthcare, there are
many mathematical and computational problems that need to be solved. 
Examples include calibration and validation of multiscale, hybrid, stochastic computational
models, forecasting algorithms, and optimal control methods. This talk will describe some
of these problems and outline a mathematical research program for the field.

\item 
\HA{\large Roland W\"uchner} (Technical University Braunschweig, Germany)

\smallskip
\noindent
{\bf Title.} Digital Twins in the Built Environment -- \underline{Challenges and Potentials}

\smallskip
\noindent
{\bf Abstract.} 
The built environment includes structures and systems with immense social, economic and environmental impacts, such as bridges, water networks, traffic and crowd management, all kinds of buildings, and conventional and renewable energy plants. Given the large amount of building material these structures require, it is of utmost importance to enable efficient material use through optimal design, long-term smart operation, and possibly retrofitting or an adapted utilization strategy. Modern numerical methods and simulations in combination with robust sensor systems can be combined to form digital twins in order to achieve this. The following special challenges of these structures must be considered and solved:

Most structures such as bridges and tall buildings are unique, i.e. built only once, which places very high demands on the predictive capabilities of numerical simulations in the design phase, but also for the evaluation of ``what-if" scenarios during operation. Hence, verification and validation of the numerical approaches are extremely important, and purely data-driven approaches for predictions are not feasible due to the lack of training data for the structures under consideration with their emerging changes during operation. Since physics-based models and simulations are indispensable and additional sensor data is available, the combination of both in hybrid FEM-ROM or FEM-ML approaches appears promising.

The actual simulation requirements – which are very immense - arise from the special characteristics of civil engineering structures. These are large, flexible structures with complex geometry and heterogeneous materials such as reinforced concrete which are immersed in their environment. Thus, the simulation techniques within the specific digital twin may require multi-scale (e.g. local concrete damage in large structures) and multi-physics (e.g. fluid-structure interaction with wind or water flows) approaches. Moreover, significant variations in material properties and boundary conditions (e.g. foundation conditions) may occur during installation and operation, which can necessitate the quantification of uncertainties.

Another aspect is the comparatively long service life. Therefore, the digital twin must be effectively adaptable even in the event of component replacement or rebuilding scenarios. This could be summarized as the need for a "maintainable digital twin", which has direct consequences on data management, software planning, interfaces between components, and a modular model setup and suitable update algorithms. Another aspect are the very different time scales (e.g. vibrations of structures versus aging of materials). While this may require suitable coupling approaches for (possibly barely) coupled phenomena with inherently very diverse time scales, the requirement for real-time calculations is much more relaxed for long-term phenomena.

The assessment of the structural condition and the detection of emerging, possibly local damage may require a digital twin model with many (preferably spatially distributed) parameters. This leads to inverse problems with very high dimensionality, which require efficient numerical methods for twinning (e.g. adjoint approaches). In cases where structures are also affected by complex multiphysics phenomena, the operation, identification and updating of the specific digital twin requires efficient approaches for the solution of coupled forward and inverse problems, which encourages the further development of coupled adjoint sensitivity analysis.
A large number of scenarios in which digital twins could be advantageously deployed and used in civil engineering can be envisioned. A non-exhaustive list includes:
\begin{itemize}
	\item Detection of aging in steel and reinforced concrete structures (bridges, high-rises, stadiums, power plants, etc.);
	\item Detection of soil/foundation softening for structures;
	
	\item Maintenance of structures and buildings;
	
	\item Condition assessment of used components for reuse scenarios to support circularity;
	
	\item Adaptation of structures to climate change and/or varying load conditions.
\end{itemize}

\item 
\HA{\large Rainald L\"ohner} (George Mason University)

\smallskip
\noindent
{\bf Title.} Digital Twins in the Built Environment -- \underline{Methods and Requirements}

\smallskip
\noindent
{\bf Abstract.} 
The term `digital twin' can mean different things in different realms - even in a specialized field like civil engineering. The data stored in the digital twin can encompass vastly different information. In some instances a BIM enumeration of parts (e.g. number of buildings in a neighbourhood) is already considered by some a digital twin. This data can be enriched by accounting for the parts (number, type, dimensions, manufacturer, … of: walls, doors, flooring, piping, HVAC systems, lighting, kitchen, etc.). And enriched further by the CAD data for visualization or product life cycle management. And enriched further by the detailed or approximate (ROM, POD, ML) computational models and sensors that monitor aspects of the built environment throughout its life. These `information rich' digital twins have emerged as a result of three megatrends: a) the pervasive use of CAD systems, b) the widespread availability and use of computational tools to `pre-compute, only then build' and `pre-compute, then operate'; and c) the emergence of precise, rugged, connected and cheap cameras and sensors that may be used to monitor structures.

The simple `data storage’ digital twins which are often used for product life management pose interesting mathematical requirements in the field of massive data storage and retrieval/comparison, pattern recognition, labeled data update, and data maintenance. 

The `information rich' of ‘high fidelity’ digital twins pose mathematical and software requirements along the complete simulation/digital twin chain:

{\bf Pre-processing:} For large, unique engineering structures, extensive predictive numerical models are created in the design phase before the actual realization of the specific product. In order to simplify the workflow to deploy a digital twin, the CAD or BIM models used for computational mechanics should be built with `digital twin readiness’ in mind. The goal is to achieve a level of automation such that if a part of the structure is changed, all metier-specific inputs (grids, material parameters, boundary conditions, ROMs, ...)  are automatically updated.

{\bf Workflow and Training:} As could be seen from above, new methods and software developments are required for the creation of integrated workflows in relation to the definition, setup, operation, updating and maintenance of digital twins for engineering structures. This also places particular demands on the work of standardization committees. In addition, (university) education and professional training must reflect the overall increase in complexity due to the growing number of models involved (possibly of various fidelity) and their dependencies and interactions in digital twins.

{\bf Solvers:} The combination of computational tools and sensors has opened the possibility to detect aging, weakening or operational hazards in real time, simplifying maintenance and improving operational efficiencies and security. In order to make this possibility a reality, a number of mathematical challenges need to be addressed.  Solvers need to be improved further (and, yes, this includes linear algebra running efficiently on advanced hardware, better preconditioners, robust nonlinear convergence and adaptive refinement). As many of the sensor/model adaptation/update tasks can be cast as optimization problems, optimization methods that operate robustly and reliably with a very high number of parameters and non-convex optimality criteria and input-based parameter adaptation need to be developed. The same applies to UQ and stochastic analysis.

{\bf Post-processing:} Given that `high fidelity’ digital twins may yield too much (or confusing) information to operators, the data must be condensed and visualized in an intuitive, metier-specific way. This will require adaptations on a case-by-case basis, but general principles based on sound mathematical foundations may help.

The talks will focus on these aspects and highlight them on two applications: detection of weakening in structures and digital twin enhanced command and control centers for crowds (DTEC4).

\item 
\HA{\large Akil Narayan} (The University of Utah)

\smallskip
\noindent
{\bf Title.} Model selection, combination, and management: Sowing with exploration and reaping with exploitation

\smallskip
\noindent
{\bf Abstract.} 
Modern simulation-based scientific models are complex and multi-faceted, involving computationally demanding physics-based modeling and discretization, reliable and robust data assimilation, and an accurate accounting for uncertainty in the face of unknown model parameter values and/or genuine stochasticity. To meet such demands, many simulations of real-world systems often involve combinations of homogenized, microscale, or reduced order model components that target accuracy of specific system subcomponents. To complicate matters, individual subcomponents may have several competing models whose query cost and accuracy payoff for system-wide prediction is opaque. For example, there may be numerous ways to identify macroscopic closure terms that are informed through analysis of a suite of microscale models. 

We discuss how computational exploration-exploitation meta-algorithms can accomplish model selection and combination for digital-twin scale applications: An exploration phase is devoted to learning about model relationships and interactions, followed by an exploitation phase that uses information learned in exploration to make decisions about optimal model selection. This framework leads to flexible procedures, capable of managing disparate models and multi-modal data in adaptive and real-time scenarios. We will identify existing theoretical guarantees for such procedures along with promising and impactful directions for new analysis and algorithm development.

\item 
\HA{\large Arvind Saibaba} (NC State)

\smallskip
\noindent
{\bf Title.} Randomized Low-rank Approximations in Scientific Computing

\smallskip
\noindent
{\bf Abstract.} 
Randomized numerical linear algebra (RandNLA) is emerging as a powerful tool to tackle large-scale linear algebra problems---such as low-rank approximations, least squares, linear systems, and matrix functions---in scientific computing and data science. RandNLA can overcome computational bottlenecks associated with classical NLA algorithms by exploiting problem structure, making judicious use of access to the underlying matrices, and enabling high-performance computing. In this talk, we will show how randomization for low-rank approximations can be used to obtain orders-of-magnitude reductions in applications of relevance to scientific computing such as sensitivity analysis, tensor decompositions, Bayesian inverse problems, and reduced order modeling.

\item 
\HA{\large Karen Veroy-Grepl} (Eindhoven University of Technology, Netherlands)

\smallskip
\noindent
{\bf Title.} Parametric Model Order Reduction in the Multi-scale Materials Setting

\smallskip
\noindent
{\bf Abstract.} 
Digital twins for multi-scale materials are extremely challenging due to the huge computational cost required for multi-scale simulations.  Such simulations are crucial in understanding structure-property relations (i.e., the effect of the microstructure on a component’s macroscopic properties), and are thus essential in the optimal design of materials for specific applications, in production process control, or to enable (for example) estimation of microstructural parameters through macroscale measurements. However, even just two-scale simulations are very computationally expensive particularly in multi-query contexts such as optimization, material design, and inverse problems. To make such analyses amenable, the microscopic simulations can be replaced by inexpensive parametric surrogate models.
 
In this talk, we present some recent work on parametric model order reduction methods to construct efficient surrogates for parametrized microstructural problems. The proposed methodology can account for parameters involving loading, material property, and geometry. The methods are tested on several (composite) microstructures, where large deformations, rotations, and finite-strain plasticity in combination with relatively large variations in the shape of inclusions are considered. The method is applied to two-scale examples, where the surrogate model achieves a high accuracy and significant speed up, thus demonstrating its potential in realistic two-scale engineering applications and material design.

\item 
\HA{\large Bart van Bloemen Waanders} (Sandia National Labs)

\smallskip
\noindent
{\bf Title.} Post-optimality sensitivities with respect to model-discrepancy for optimization of complex systems

\smallskip
\noindent
{\bf Abstract.} 
Digital twins require a fusion of numerical modeling and outer-loop analysis to support decision-making for complex systems.  In this talk,  PDE-constrained optimization serves as the outer-loop analysis and the goal is to achieve an optimization solution at the highest possible model fidelity.  However, implementing large scale optimization (objective functions, adjoints, Hessians) is often intractable for high-fidelity models.  To that end, we present a novel use of post-optimality sensitivities with respect to model discrepancy using lower-fidelity PDEs to achieve a higher-fidelity optimization solution.  We demonstrate our approach on various numerical exemplars, including control of Stokes to emulate control of Navier Stokes.

\item 
\HA{\large Carol Woodward} (Lawerence Livermore National Lab)

\smallskip
\noindent
{\bf Title.} Developing and Deploying Software for Use in Multi-Temporal-Scale Applications on High Performance Computing Systems

\smallskip
\noindent
{\bf Abstract.} 
In collaboration with Cody J. Balos, David J. Gardner, Steven B. Roberts, and Daniel R. Reynolds

Effective digital twins rely on robust numerical software that can be utilized in multiple computing environments. This software must provide implementations of numerous methods including algorithms for optimization, algebraic solvers, time evolution, and data assimilation, among others.  The SUNDIALS software library provides efficient, adaptive time integration algorithms for use in application codes with time dynamics.  Through the Exascale Computing Project, the team has rearchitected the software platform to provide new flexibilities for accommodating use on hybrid computer architectures, use in varied application contexts, and use of several algebraic solver and data structure classes.  Through the DOE SciDAC program the team has researched new methods for multiscale applications as well as deployed the software to several applications.
In this talk I will discuss considerations for development of software for public distribution and use on high performance computing systems using examples from SUNDIALS’ experience.  The importance of the three areas of method research, capability development, and applications deployment will be discussed along with recent results and application examples.
This work was performed under the auspices of the U.S. Department of Energy by Lawrence Livermore National Laboratory under Contract DE-AC52-07NA27344. Lawrence Livermore National Security, LLC. LLNL-ABS-854222.

\item 
\HA{\large Enrique Zuazua} (Friedrich-Alexander-Universit\"at Erlangen-Nürnberg, Germany)

\smallskip
\noindent
{\bf Title.} Control and Machine Learning

\smallskip
\noindent
{\bf Abstract.} 
Control theory and Machine Learning share common objectives, as evident in Norbert Wiener's definition of ``Cybernetics" as ``The science of control and communication in animals and machines."

The synergy between these fields is reciprocal. Control theory tools enhance our comprehension of the efficacy of certain Machine Learning algorithms and offer insights  for their enhancement. However, this often bounces  intricate queries back.

Consider the control of a linear finite-dimensional system—as example. A sharp mathematical solution exists: it suffices to ensure the Kalman matrix's rank is full. Yet, computing these matrices in high dimensions presents a new challenge. DeepMind has made remarkable contributions, introducing artificial intelligence solutions to the old problem of matrix multiplication.

The interplay between Control and Machine Learning opens up a new captivating  scientific  landscape to be explored but this can become a labyrinthine task. And this  is surely part of the overall ambitious program of developing Digital Twins technologies.
In this talk, we will present some of the contributions from our team at the interface between Control  and Machine Learning, which can modestly contribute to this noble and complex task.

\end{enumerate}

\newpage
\begin{center}
 \underline{\bf\Large Applications Panel}
\end{center}

\vspace{1cm}

\begin{enumerate}[1.]
\setlength\itemsep{2em}
\item 
\HA{\large Aldo Badano} (FDA)

\smallskip
\noindent{\bf Title:} Digital Twins and Other Synthetic Humans in the Regulatory Evaluation of Medical Devices 

\smallskip
\noindent{\bf Abstract:} Applications of digital twin technology in regulatory science differ from applications in other industries because of data complexity, data sparsity, data privacy, and data ownership issues. To realize the potential of digital twins in biomedicine, several scientific challenges need to be addressed including commutability, bias mitigation, wide-sense generalizability, and platforms, standards, and protocols for in silico databanks.

\item 
\HA{\large Stefan Boschert} (SIEMENS)

\smallskip
\noindent{\bf Title:} Challenges for Digital Twins in Industrial Applications

\smallskip
\noindent{\bf Abstract:}
Digital Twins have become very popular in the last decade. Its development strongly benefited from advances in simulation technology and the rising availability of data so many products are already already equipped with Digital Twin technology. However, there are several challenges from an industry point of view, that need to be tackled before Digital Twins can leverage their full potential. In the talk I will give some examples of potential industrial applications, where existing technologies come to its limits in actual realization.

\item 
\HA{\large Brent Bartlett} (NVIDIA)

\smallskip
\noindent{\bf Title:} Creating Geospatial Workflows with Digital Twin

\smallskip
\noindent{\bf Abstract:}
Edge sensors are quickly becoming indispensable to many different use cases and can provide rich information streams to enrich virtual representations of the real world. The amount of data these sensors create is growing at an incredible rate and new system architectures are needed to effectively update Digital Twins. This talk will present how NVIDIA Holoscan can be utilized to create geospatial processing workflows that link sensor data to a digital twin with examples of FMV and SAR data streaming into and NVIDIA Omniverse digital twin.

\item 
\HA{\large Juan R. Cebral} (George Mason University)

\smallskip
\noindent{\bf Title:} Biomedical Digital Twins from Virtual Patient Populations

\smallskip
\noindent{\bf Abstract:} Personalization of surgical procedures and medical therapies requires understanding the most likely responses of each individual patient to the different treatment options. However, assessing treatment outcomes and potential complications depends on the patient’s individual anatomical and physiologic characteristics, which are difficult to determine with certainty. Building a digital twin of a given patient that reproduces the salient features associated with therapeutic outcomes provides the opportunity to make more informed, objective, and personalized clinical decisions. 
To achieve this goal, we envision building a virtual patient population (i.e. computational models) with varying anatomical and physiologic characteristics and pre-compute synthetic images and/or physiologic curves that could be directly compared to actual patient data and measurements. Next, the virtual patient with the closest (i.e. most similar) imaging and physiologic data can be identified using statistical and/or machine learning approaches. Once so identified, this (approximate) digital twin can be used to assess different therapeutic options and inform the clinical decision-making process.

As an example, we investigate different interventional procedures used for occlusive stroke treatment. The outcomes of these procedures largely depend on the patient's individual brain vascular collateralization, which provide alternative avenues of flow to the area of the brain affected by an arterial occlusion. The human brain contains collateral vessels (i.e. communicating arteries) that join the main feeding arteries at the level of the circle of Willis at the base of the brain, as well as a network of smaller collaterals at the surface of the brain called pial collaterals or leptomeningeal arteries. However, there is a large variability of these collateral vessels among the population. Even tough medical images such as computed tomography angiography (CTA) can depict the larger collaterals, their resolution is limited, and smaller vessels are not clearly visualized. Nevertheless, temporal imaging such as multi-phase CTA provide valuable information about the delayed filling of smaller arteries feeding the affected area of the brain (through collaterals). This information can be directly used to identify a digital twin within a virtual patient population with varying levels of brain collateralization that best matches (i.e. explains) the patient’s temporal curves. Finally, the patient data can be augmented with data from the digital twin and/or different treatment options could be investigated using the digital twin. Thus, improved understanding of the underlying anatomy and effects of interventions on a patient-specific basis using the corresponding digital twin has the potential to further improve patient management.

The proposed strategy poses a number of mathematical and computational challenges, including: 1) creating virtual patient populations with the correct anatomical variability (e.g. using constrained constructive optimization methods to extend vascular models), 2) creating synthetic data that can be compared to patient data (e.g. performing transport simulations to create virtual angiograms and quantifying specific features from medical images and virtual images), 3) identifying digital twins that best match patients’ data (e.g. using statistical correlations, machine learning models, etc.), 4) using the digital twins to assess best options for a patient (e.g. performing numerical simulations of medical procedures, identifying statistical associations between anatomical and physiological characteristics and intervention outcomes, creating machine learning models to predict outcomes, etc.)

\item \HA{\large Sharon Di} (Columbia University)

\smallskip
\noindent{\bf Title:}
AI for Urban Transportation Digital Twin

\smallskip
\noindent{\bf Abstract:}
Transportation digital twins have become increasingly popular tools to improve traffic efficiency and safety. However, the majority of effort nowadays is focused on the “eyes” of the digital twin, which is object detection using computer vision. I believe the key to empowering the intelligence of a transportation digital twin lies in its “brain,” namely, how to utilize the information extracted from various sensors to infer traffic dynamics evolution and devise optimal  control and management strategies with real-time feedback to guide the transportation ecosystem toward a social optimum. 

My research aims to employ tools including machine learning and game theory to develop an urban transportation digital twin, leveraging data collected from the NSF PAWR COSMOS city-scale wireless testbed being deployed in West Harlem next to the Columbia campus. In this talk, I will primarily focus on two solutions: (1) scientific machine learning that leverages both domain knowledge and available data, and (2) mean field game that bridges the gap between micro- and macroscopic behaviors of multi-agent dynamical systems. In the first topic, physics-informed deep learning will be introduced and applied to traffic state estimation and uncertainty quantification. In the second topic, I will introduce how to model behaviors of new actors (e.g., a large number of autonomous vehicles) in a transportation system and their interaction with existing actors (e.g., human-driven vehicles).

\item 
\HA{\large Mustafa Megahed} (ESI Group)

\smallskip
\noindent{\bf Title:} Hybrid Twins for divers Data Sets in Industrial Applications

\smallskip
\noindent{\bf Abstract:}
Given the very large differences in data sources, acquisition frequency, accuracy and formats, the large variety in reasons behind using digital twins and the large assortment of data analysis and machine learning tools available the question arises ``How straight forward is it to apply digital twins in industrial applications?". This talk will explore these aspects and how they impact industrial implementation. The hybrid twin concept is introduced as applied to material behavior, manufacturing outcome and vehicle safety predictions.

\end{enumerate}

\newpage
\begin{center}
 \underline{\bf\Large Contributed Posters and Demos}
\end{center}

\vspace{1cm}
\begin{enumerate}
\item \HA{\large Facundo Airaudo} (George Mason University)\\
\textit{On the use of risk measures in digital twins to identify weaknesses in structures}
\item \HA{\large Reza Akbarian Bafghi} (University of Colorado Boulder)\\
\textit{Accelerating PINNs implementations in PyTorch and TensorFlow}
\item \HA{\large Antonio Álvarez López} (Universidad Autónoma De Madrid)\\
\textit{Optimized classification with neural ODEs}
\item \HA{\large Mohammad Atif} (Brookhaven National Laboratory)\\
\textit{Towards enabling digital twins capabilities for a cloud chamber}
\item \HA{\large Charles Beall} (Stevens Institute of Technology)\\
\textit{Randomized reduced basis methods for advection-diffusion problems}
\item \HA{\large Prajakta Bedekar} (Johns Hopkins University and National Institute of Standards and Technology)\\
\textit{Time-dependent antibody kinetics for previously infected and vaccinated individuals via graph-theoretic modeling}
\item \HA{\large Filip Belik} (University of Utah)\\
\textit{Dynamic arterial bulk conductivity from blood pressure}
\item \HA{\large Md Rezwan Bin Mizan}	(University of Houston)\\
\textit{Reduced order modeling of the Kuramoto-Sivashinsky equation using proper orthogonal decomposition}
\item \HA{\large Mingchao Cai} (Morgan State University)\\
\textit{A combination of physics-informed neural networks with the fixed-stress splitting iteration for solving Biot's model} 
\item \HA{\large Sean Carney} (George Mason University)\\
\textit{PDE constrained optimization with distributional uncertainty}
\item \HA{\large Rujeko Chinomona} (Temple University)\\
\textit{Real-time interactive simulation of reality with neuroscience applications}
\item \HA{\large Kamala Dadashova} (North Carolina State University)\\
\textit{Parameter subset selection for identifiability analysis in mPBPK model}
\item \HA{\large Alejandro Diaz} (Rice University)\\
\textit{Nonlinear manifold reduced order models with domain decomposition}
\item \HA{\large Emmanuel Fleurantin} (George Mason University and University of North Carolina at Chapel Hill)\\
\textit{A data driven study of the drivers of stratospheric circulation via reduced order modeling and data assimilation}
\item \HA{\large Jithin George} (Northwestern University)\\
\textit{Walking into the complex plane to ``order" better time integrators}
\item \HA{\large Sebastian Gutierrez Hernandez} (Georgia Institute of Technology)\\
\textit{Parametric geodesics and control problems in Wasserstein space}
\item \HA{\large Sebastian Kaltenbach} (Harvard University)\\
\textit{Korali - part 2: Bayesian inference}
\item \HA{\large Kai Fung (Kelvin) Kan} (Emory University)\\
\textit{LSEMINK: a modified Newton-Krylov method}
\item \HA{\large Prashant Kambali} (Villanova University)\\
\textit{Hybrid Modeling Approaches for Dynamic Systems: Bridging Nonlinear Physics, Machine Learning, and Expert Insights in Engineering and Medicine }
\item \HA{\large Rohit Khandelwal} (George Mason University)\\
\textit{A discontinuous Galerkin method for optimal control of the obstacle problem}
\item \HA{\large Hyunah Lim} (University of Maryland College Park)\\
\textit{Mathematical modelling of the impacts of screening and vaccination in HPV}
\item \HA{\large Anton Malandii} (Stony Brook University)\\
\textit{Adaptive importance sampling for topology optimization in uncertain environment}
\item \HA{\large Zhaobin Mo} (Columbia University)\\
\textit{Physics-Informed Deep Learning for Traffic State Estimation: A Survey and the Outlook}
\item \HA{\large Soyoung Park} (University of Maryland)\\
\textit{Mathematical modelling of the impacts of screening and vaccination in HPV}
\item \HA{\large Aseem Milind Pradhan} (George Mason University)\\
\textit{A digital cerebral vasculature database for stroke management}
\item \HA{\large Rhea Shroff} (University of Florida)\\
\textit{Accelerating the computation of tensor $z$-eigenvalues}
\item \HA{\large Priscila Silva} (University of Massachusetts Dartmouth)\\
\textit{Predictive system resilience modeling}
\item \HA{\large Moyi Tian} (Brown University)\\
\textit{Community robustness under edge addition in synthetic and empirical temporal networks}
\item \HA{\large Shanyin Tong} (Columbia University)\\
\textit{Large deviation theory-based adaptive importance sampling for rare events in high dimensions}
\item \HA{\large Deepanshu Verma} (Emory University)\\
\textit{Neural Network Approaches for Parameterized Optimal Control}
\end{enumerate}

\newpage
\begin{center}
 \underline{\bf\Large Poster Titles and Abstracts}
\end{center}
\vspace{1cm}
\begin{center}
\vspace{0.2cm}
{\bf\large On the use of risk measures in digital twins to identify weaknesses in structures}
\vspace{0.2cm}\\
{\bf\large \underline{Facundo Airaudo}$^1$}  

\vspace{0.5cm}
{\large $^1$George Mason University
}
\end{center}

\vspace{0.5cm}
Given measurements from sensors and a set of standard forces, an optimization based approach to identify weakness in structures is introduced. The key novelty lies in letting the load and measurements to be random variables. Subsequently the conditional-value-at-risk (CVaR) is minimized subject to the elasticity equations as constraints. CVaR is a risk measure leads to minimization of rare and low probability events which the standard expectation cannot. The optimization variable is the (deterministic) strength factor which appears as a coefficient in the elasticity equation, thus making the problem nonconvex. Due to uncertainty, the problem is high dimensional and, due to CVaR, the problem is nonsmooth. An adjoint based approach is developed with quadrature in the random variables. Numerical results are presented in the context of a weakened region in a plate and a large structure with trusses similar to those used in solar arrays or cranes.


\vspace{1cm}
\hrulefill
\vspace{1cm}
\begin{center}
\vspace{0.2cm}
{\bf\large Accelerating PINNs implementations in PyTorch and TensorFlow}
\vspace{0.2cm}\\
{\bf\large Reza Akbarian Bafghi$^1$}  

\vspace{0.5cm}
{\large $^1$University of Colorado Boulder
}
\end{center}

\vspace{0.5cm}
We want to introduce two Python packages that accelerate the implementation of Physics-Informed Neural Networks (PINNs) using the PyTorch and TensorFlow frameworks. These packages streamline user interaction by abstracting PDE issues.


\vspace{1cm}
\hrulefill

\vspace{1cm}
\begin{center}
\vspace{0.2cm}
{\bf\large Optimized classification with neural ODEs}
\vspace{0.2cm}\\
{\bf\large Antonio Álvarez López$^1$}  

\vspace{0.5cm}
{\large $^1$Universidad Autónoma De Madrid
}
\end{center}

\vspace{0.5cm}
Classification of $O(N)$ points becomes a simultaneous control problem when viewed through the lens of neural ordinary differential equations (neural ODEs), the time-continuous limit of residual networks. We estimate the likelihood that each number of neurons below $O(N)$ may be required for that task in the worst-case scenario, when the points are independently and uniformly distributed in $[0,1]^d$. Notably, we quantify how this probability evolves with respect to $N$ when the dimension $d$ grows. Furthermore, assuming a generic condition on the point distribution, we propose a new constructive algorithm for classification that steers entire clusters of $d$ points at once, leading to a reduced maximal complexity of $O(N/d)$ neurons.


\vspace{1cm}
\hrulefill

\vspace{1cm}
\begin{center}
\vspace{0.2cm}
{\bf\large Towards enabling digital twins capabilities for a cloud chamber}
\vspace{0.2cm}\\
{\bf\large Mohammad Atif$^1$}  

\vspace{0.5cm}
{\large $^1$Brookhaven National Laboratory
}
\end{center}

\vspace{0.5cm}
Particle-resolved direct numerical simulations (PR-DNS), which resolve not only the smallest turbulent eddies but also track the development and motion of individual particles, are arguably an essential tool for exploring aerosol-cloud-turbulence interactions at the fundamental level. For instance, PR-DNS may complement experimental facilities designed to study key physical processes in a controlled environment and therefore serve as digital twins for such cloud chambers. In this poster we present our ongoing work aimed at enabling the use of a PR-DNS model for this purpose. We consider two approaches: traditional HPC techniques and emerging machine learning methods. Future research directions are outlined as well.


\vspace{1cm}
\hrulefill

\vspace{1cm}
\begin{center}
\vspace{0.2cm}
{\bf\large Randomized reduced basis methods for advection-diffusion problems}
\vspace{0.2cm}\\
{\bf\large Charles Beall$^1$}  

\vspace{0.5cm}
{\large $^1$Stevens Institute of Technology
}
\end{center}

\vspace{0.5cm}
In this project, we develop randomized model reduction methods for advection-diffusion problems with sharply discontinuous source terms. To study such problems, we must solve the advection-diffusion equation, a partial differential equation (PDE) used to model systems such as a liquid dye being dissolved in a flowing fluid, or the combination of heat conduction and convection through a medium. This PDE arises often in the sciences, from fluid dynamics [1] to semiconductor physics [2], although in such contexts, the equation is almost always unsolvable by hand, so we rely on computer algorithms to efficiently obtain accurate approximations to the solutions. Our goal is to obtain a faster result than with direct numerical simulations like finite difference or finite element schemes. We employ randomized methods from data science to allow for parallel-in-time computation [3] and generation of a reduced order model. Compared to direct simulations, the reduced model produces a solution space of much lower dimension, meaning the computational complexity is greatly reduced. Thus, simulations with the reduced model allow for a further speedup in computational runtime, with the added benefit of maintaining accuracy. As a novel contribution, we consider the case of sharp discontinuities in source functions, partition the time domain into overlapping subintervals with overlap around the discontinuities, and construct a reduced basis on each subinterval. This allows for the construction of a reduced solution that combines information from the reduced bases, rather than relying on a single reduced basis to capture information throughout the time domain. We present a test case to show that this approach can provide significant improvements in accuracy compared to the construction of only one reduced basis over the whole time domain. \\

\textbf{References:}
\begin{enumerate}
\item[1] Bird, R. B., Stewart, W. E. \& Lightfoot, E. N. (2007). Transport Phenomena (Revised second). John Wiley \& Sons. 
\item[2] New, O. (2004) Derivation and numerical approximation of the quantum drift diffusion model for semiconductors, Jour. Myan. Acini. Arts \& Sc., Vol. II (Part Two), No. 5. 
\item[3] Schleuß, J., Smetana, K., \& ter Maat, L. (2022). Randomized quasi-optimal local approximation spaces in time. arXiv preprint arXiv:2203.06276. To appear in SIAM J. Sci. Comput., 2023.
\end{enumerate}

\vspace{1cm}
\hrulefill

\vspace{1cm}
\begin{center}
\vspace{0.2cm}
{\bf\large Time-dependent antibody kinetics for previously infected and vaccinated individuals via graph-theoretic modeling}
\vspace{0.2cm}\\
{\bf\large Prajakta Bedekar$^1$}  

\vspace{0.5cm}
{\large $^1$Johns Hopkins University and National Institute of Standards and Technology
}
\end{center}

\vspace{0.5cm}
Modeling the deterioration of antibody levels is paramount to understanding the time-dependent viral response to infections, vaccinations, or a combination of the two. These events have been studied experimentally, but also benefit from a rigorous mathematical underpinning. Disease/vaccination prevalence in the population and time-dependence on a personal scale simultaneously affect antibody levels, interact non-trivially, and pose considerable modeling challenges. We propose a time-inhomogeneous Markov chain model for event-to-event transitions coupled with a probabilistic framework for post-infection or post-vaccination antibody kinetics. This approach is ideal to model sequences of infections and vaccinations, or personal trajectories in a population. We use synthetic data to demonstrate the modeling process as well as estimation of transition probabilities. This work is an important step towards a comprehensive understanding of antibody kinetics that will allow us to simulate and predict real-world disease response scenarios. Modeling the deterioration of antibody levels is paramount to understanding the time-dependent viral response to infections, vaccinations, or a combination of the two. These events have been studied experimentally, but also benefit from a rigorous mathematical underpinning. Disease/vaccination prevalence in the population and time-dependence on a personal scale simultaneously affect antibody levels, interact non-trivially, and pose considerable modeling challenges. We propose a time-inhomogeneous Markov chain model for event-to-event transitions coupled with a probabilistic framework for post-infection or post-vaccination antibody kinetics. This approach is ideal to model sequences of infections and vaccinations, or personal trajectories in a population. We use synthetic data to demonstrate the modeling process as well as estimation of transition probabilities. This work is an important step towards a comprehensive understanding of antibody kinetics that will allow us to simulate and predict real-world disease response scenarios. 

\vspace{1cm}
\hrulefill

\vspace{1cm}
\begin{center}
\vspace{0.2cm}
{\bf\large Dynamic arterial bulk conductivity from blood pressure}
\vspace{0.2cm}\\
{\bf\large Filip Belik$^1$}  

\vspace{0.5cm}
{\large $^1$University of Utah
}
\end{center}

\vspace{0.5cm}
Impedance cardiography is a non-invasive procedure in which an alternating electrical current is passed into the human body to measure total conductivity. Currents are primarily conducted through the arteries and blood as they generally exhibit lower resistivity compared to the muscle, skin, and fat. This conductivity can then be used to study important physiological parameters such as stroke volume, heart rate, and cardiac output. We propose an analytical forward model to understand the complicated relationship between blood pressure measured at the upper arm and the electrical conductivity of the radial artery at the wrist. This model involves several steps. First, the transport of a pressure wave down the brachial artery, through the bifurcation of the brachial artery into the ulnar and radial arteries, and then the transport down the radial artery to the wrist. Next, reduced Navier Stokes equations are matched with a moving boundary to model the pulsatile nature of blood flow alongside the elasticity of the tube wall. Finally, Maxwell-Fricke equations are used alongside the shear stress-induced orientation and deformation of red blood cells to calculate the bulk conductivity of the blood.


\vspace{1cm}
\hrulefill

\vspace{1cm}
\begin{center}
\vspace{0.2cm}
{\bf\large Reduced order modeling of the Kuramoto-Sivashinsky equation using proper orthogonal decomposition}
\vspace{0.2cm}\\
{\bf\large \underline{Md Rezwan Bin Mizan}$^1$, Ilya Timofeyev$^1$, Maxim Olshansky$^1$, and Alexander Mamonov$^1$}  

\vspace{0.5cm}
{\large $^1$University of Houston
}
\end{center}

\vspace{0.5cm}
Proper Orthogonal Decomposition (POD) techniques are widely utilized in various scientific fields to streamline complex, data-driven spatio-temporal dynamics. In this study, we present a POD-Reduced Order Model (POD-ROM) for the one-dimensional Kuramoto-Sivashinsky (KS) equation, known for its chaotic attractor. Our FOM simulation was conducted over an extended time period to capture snapshots within this attractor. We constructed our basis using truncated Singular Value Decomposition (SVD), effectively representing the attractor. This single basis was then employed to construct our ROM using Galerkin projection. We demonstrate that the POD- ROM can accurately reproduce the statistical features of the KS equation’s attractor, independent of initial conditions. Additionally, we show that this basis allows short-term predictions of the KS dynamics. Our findings suggest that a POD-ROM, based on a single basis, is effective for both long-term and short-term predictions of the chaotic KS dynamics. Furthermore, we have developed a criteria for the automatic selection of the reduced model’s dimension, based on the cumulative variance of the singular values.

Keywords: Kuramoto-Sivashinsky Equation, Proper Orthogonal Decomposition (POD), spatio-temporal dynamics, Reduced Order Model (ROM)

\vspace{1cm}
\hrulefill

\vspace{1cm}
\begin{center}
\vspace{0.2cm}
{\bf\large A combination of physics-informed neural networks with the fixed-stress splitting iteration for solving Biot's model}
\vspace{0.2cm}\\
{\bf\large Mingchao Cai$^1$}  

\vspace{0.5cm}
{\large $^1$Morgan State University
}
\end{center}


\vspace{1cm}
\hrulefill

\vspace{1cm}
\begin{center}
\vspace{0.2cm}
{\bf\large PDE constrained optimization with distributional uncertainty}
\vspace{0.2cm}\\
{\bf\large Sean Carney$^1$}  

\vspace{0.5cm}
{\large $^1$George Mason University
}
\end{center}

\vspace{0.5cm}
A critical task when seeking to quantify and minimize the risk of rare, outlier events, such as structural failures or climate catastrophes, is to account for uncertainty in the underlying physical models. Uncertainty in physical systems is modeled with random variables, however, in practice there may be some nontrivial ambiguity in the underlying probability distribution from which they are sampled. This work describes an analytic framework and computational methods that guarantee minimizers are robust to such ambiguities.


\vspace{1cm}
\hrulefill

\vspace{1cm}
\begin{center}
\vspace{0.2cm}
{\bf\large Real-time interactive simulation of reality with neuroscience applications}
\vspace{0.2cm}\\
{\bf\large Rujeko Chinomona$^1$}  

\vspace{0.5cm}
{\large $^1$Temple University
}
\end{center}

\vspace{0.5cm}
In this project, we delve into some crucial aspects of digital twins, focusing on real-time interaction and the simulation of real-time data. Unlike traditional high performance computing, which typically involves running simulations and creating visualizations afterward, our Virtual Simulation Of Reality (VISOR) project enables real-time interaction and immediate visualization, fostering an intuitive understanding of the physical system and allowing on-the-fly adjustments. As part of the VISOR project, we have developed a computational neuroscience tool called Neuroscientific Virtual Simulation Of Reality (Neuro-VISOR). This tool simulates the Hodgkin-Huxley equations, which describe the active electrical signal processing in neurons. Our primary goal is to create real-time, interactive simulations tailored to assist researchers in their neuroscientific experiments, ultimately working towards the development of digital twins for neuroscience research. We also address some of the numerical challenges related to computing accurate solutions to the Hodgkin-Huxley equations in real-time, ensuring a seamless and fluid visualization experience. During our demonstration, participants will have the opportunity to explore Neuro-VISOR. Equipped with a virtual reality (VR) headset, they can visualize the current system state in real-time and make instant adjustments. Some noteworthy features of the current simulation model include the capability to directly modify voltage at specific neuron points, insert a clamp on a designated neuron point, and create synaptic connections between different neurons.

\vspace{1cm}
\hrulefill

\vspace{1cm}
\begin{center}
\vspace{0.2cm}
{\bf\large Parameter subset selection for identifiability analysis in mPBPK model}
\vspace{0.2cm}\\
{\bf\large Kamala Dadashova$^1$}  

\vspace{0.5cm}
{\large $^1$North Carolina State University
}
\end{center}

\vspace{0.5cm}
We consider a minimal PBPK (mPBPK) model of the brain for antibody therapeutics. This model is the reduced form of an previous multi-species platform brain PBPK model. The original PBPK model consists of 100 differential equations, whereas the mPBPK model contains 16 differential equations, improving the speed of simulations. The reduced model includes 31 parameters, and their values are taken from the original brain PBPK model. We implement a local sensitivity-based parameter subset selection (PSS) algorithm to determine identifiable parameters based on verified threshold values. We also compare results to a more standard PSS algorithm. One objective is to determine which parameters have the greatest influence on model predictions and how changes in these parameters affect model responses. This information can be utilized to improve the model's performance and facilitate subsequent uncertainty analysis. We systematically verify obtained results using qualitative methods and quantitative techniques based on energy statistics. 

\vspace{1cm}
\hrulefill

\vspace{1cm}
\begin{center}
\vspace{0.2cm}
{\bf\large Nonlinear manifold reduced order models with domain decomposition}
\vspace{0.2cm}\\
{\bf\large Alejandro Diaz$^1$}  

\vspace{0.5cm}
{\large $^1$Rice University
}
\end{center}

\vspace{0.5cm}
This poster discusses the integration of nonlinear-manifold reduced order models (NM-ROMs) with domain decomposition (DD). NM-ROMs approximate the full order model (FOM) state in a nonlinear-manifold by training a shallow, sparse autoencoder using FOM snapshot data. These NM-ROMs can be advantageous over linear-subspace ROMs (LS-ROMs) for problems with slowly decaying Kolmogorov n-width. However, the number of NM-ROM parameters that need to trained scales with the size of the FOM. Moreover, for "extreme-scale" problems, the storage of high-dimensional FOM snapshots alone can make ROM training expensive. To alleviate the training cost, DD is applied to the FOM, NM-ROMs are computed on each subdomain, and are coupled to obtain a global NM-ROM. This approach has several advantages: Subdomain NM-ROMs can be trained in parallel, each involve fewer parameters to be trained than global NM-ROMs, require smaller subdomain FOM dimensional training data, and subdomain NM-ROMs can be tailored to subdomain-specific features of the FOM. The shallow, sparse architecture of the autoencoder used in each subdomain NM-ROM allows application of hyper-reduction (HR), reducing the complexity caused by nonlinearity and yielding computational speedup of the NM-ROM. The proposed DD NM-ROM with HR approach is numerically compared to a DD LS-ROM with HR on 2D steady-state Burgers' equation, showing an order of magnitude improvement in accuracy of the proposed DD NM-ROM over the DD LS-ROM.


\vspace{1cm}
\hrulefill

\vspace{1cm}
\begin{center}
\vspace{0.2cm}
{\bf\large A data driven study of the drivers of stratospheric circulation via reduced order modeling and data assimilation}
\vspace{0.2cm}\\
{\bf\large Emmanuel Fleurantin$^1$}  

\vspace{0.5cm}
{\large $^1$George Mason University and University of North Carolina at Chapel Hill
}
\end{center}

\vspace{0.5cm}
Stratospheric dynamics are strongly affected by the absorption/emission of radiation in the Earth's atmosphere and waves which propagate upward from the troposphere perturbing stratospheric zonal flow. Reduced order models of stratospheric wave-zonal interactions, which parameterize radiative drivers and tropospheric perturbation, have been used in the past to study interannual variability in stratospheric zonal winds and sudden stratospheric warming (SSW) events. These models are most sensitive to two main parameters $\Lambda$, which forces the mean radiative zonal wind, and $h$ which acts as a perturbation parameter representing the effect of waves propagating upward from the troposphere. These parameters have been taken to be constant or have been given various time dependent parameterizations to account for seasonal, solar, and other possible forcings. In this work, we take one such reduced order model (Ruzmaikin 2002) in conjunction with 20 years of atmospheric reanalysis data provided by the European Center for Medium Weather Forecasting (ECMWF) to estimate $\Lambda$ and $h$ using various methods in data assimilation. We employ a particle filter to estimate constants in previous parameterizations of $\Lambda$ and $h$ as well as ensemble smoothing with multiple data assimilation (ESMDA) to estimate time series for $\Lambda(t)$ and $h(t)$ such that model output closely matches the reanalysis data. We provide some analysis of the resulting time series, compare them to previous parameterizations, and examine them around known historical SSW events. This allows a data driven examination of these important parameters through the lens and tractability of a reduced order model.


\vspace{1cm}
\hrulefill

\vspace{1cm}
\begin{center}
\vspace{0.2cm}
{\bf\large Walking into the complex plane to ``order" better time integrators}
\vspace{0.2cm}\\
{\bf\large Jithin George$^1$}  

\vspace{0.5cm}
{\large $^1$Northwestern University
}
\end{center}

\vspace{0.5cm}
Most numerical methods for time integration use real time steps. Complex time steps provide an additional degree of freedom, as we can select the magnitude of the step in both the real and imaginary directions. By time stepping along specific paths in the complex plane, integrators can gain higher orders of accuracy or achieve expanded stability regions. We show how to derive these paths for explicit and implicit methods, discuss computational costs and storage benefits, and demonstrate clear advantages for complex-valued systems like the Schrodinger equation. We also explore how complex time stepping also allows us to break the Runge-Kutta order barrier, enabling 5th order accuracy using only five function evaluations for real-valued differential equations.


\vspace{1cm}
\hrulefill

\vspace{1cm}
\begin{center}
\vspace{0.2cm}
{\bf\large Parametric geodesics and control problems in Wasserstein space}
\vspace{0.2cm}\\
{\bf\large Sebastian Gutierrez Hernandez$^1$}  

\vspace{0.5cm}
{\large $^1$Georgia Institute of Technology
}
\end{center}

\vspace{0.5cm}
In the work ``Parameterized Wasserstein Hamiltonian flow", the authors propose a parametric formulation to solve Wasserstein Hamiltonian Flows numerically. In this poster, we offer a short introduction to this topic and some recent developments. In particular, we show that when the considered distributions are Gaussians, this formulation can covert the Wasserstein geodesic and barycenter problems into equivalent problems in Euclidean space, which allows us to design efficient algorithms to compute the Wasserstein barycenter even in higher dimensions. Additionally, we apply the parametric formulation to control problems in the Wasserstein space.

\vspace{1cm}
\hrulefill

\vspace{1cm}
\begin{center}
\vspace{0.2cm}
{\bf\large Korali: Stochastic Optimization and Bayesian Inference}
\vspace{0.2cm}\\
{\bf\large Sebastian Kaltenbach$^1$}  

\vspace{0.5cm}
{\large $^1$Harvard University
}
\end{center}

\vspace{0.5cm}
During this demo session, we are presenting Korali (https://www.cse-lab.ethz.ch/korali/), a high-performance framework for uncertainty quantification, optimization, and deep reinforcement learning. Korali’s engine provides support for large-scale HPC systems and a multi-language interface compatible with distributed computational models. This presentation will focus on addressing (stochastic) optimization challenges as well as Bayesian inverse problems with Korali and showcase algorithms such as CMA-ES and Transitional Markov Chain Monte Carlo. In case evaluating the original model is expensive, we demonstrate how a surrogate model can be used instead. Moreover, we show how Korali is used within DCoMEX, a European High Performance Computing Joint Undertaking project.


\vspace{1cm}
\hrulefill

\vspace{1cm}
\begin{center}
\vspace{0.2cm}
{\bf\large LSEMINK: a modified Newton-Krylov method}
\vspace{0.2cm}\\
{\bf\large Kai Fung (Kelvin) Kan$^1$}  

\vspace{0.5cm}
{\large $^1$Emory University
}
\end{center}

\vspace{0.5cm}
This paper introduces LSEMINK, an effective modified Newton-Krylov algorithm geared toward minimizing the log-sum-exp function for a linear model. Problems of this kind arise commonly, for example, in geometric programming and multinomial logistic regression. Although the log-sum-exp function is smooth and convex, standard line search Newton-type methods can become inefficient because the quadratic approximation of the objective function can be unbounded from below. To circumvent this, LSEMINK modifies the Hessian by adding a shift in the row space of the linear model. We show that the shift renders the quadratic approximation to be bounded from below and that the overall scheme converges to a global minimizer under mild assumptions. Our convergence proof also shows that all iterates are in the row space of the linear model, which can be attractive when the model parameters do not have an intuitive meaning, as is common in machine learning. Since LSEMINK uses a Krylov subspace method to compute the search direction, it only requires matrix-vector products with the linear model, which is critical for large-scale problems. Our numerical experiments on image classification and geometric programming illustrate that LSEMINK considerably reduces the time-to-solution and increases the scalability compared to geometric programming and natural gradient descent approaches. It has significantly faster initial convergence than standard Newton-Krylov methods, which is particularly attractive in applications like machine learning. In addition, LSEMINK is more robust to ill-conditioning arising from the nonsmoothness of the problem. We share our MATLAB implementation at https://github.com/KelvinKan/LSEMINK.


\vspace{1cm}
\hrulefill

\vspace{1cm}
\begin{center}
\vspace{0.2cm}
{\bf\large Hybrid Modeling Approaches for Dynamic Systems: Bridging Nonlinear Physics, Machine Learning, and Expert Insights in Engineering and Medicine }
\vspace{0.2cm}\\
{\bf\large Prashant Kambali $^1$}  

\vspace{0.5cm}
{\large $^1$ Villanova University
}
\end{center}

\vspace{0.5cm}
This research endeavors to pioneer advanced techniques for modeling dynamic systems by synergistically incorporating nonlinear physics, machine learning methodologies, and expert insights. The interdisciplinary approach aims to enhance the accuracy and predictive capabilities of models in diverse fields such as engineering and medicine. By amalgamating the precision of nonlinear physics with the adaptability of machine learning algorithms and the nuanced understanding provided by expert insights, this study seeks to unlock new frontiers in system modeling. The outcomes hold significant promise for optimizing design processes, improving decision-making in complex systems, and advancing personalized medical diagnostics and treatments. 


\vspace{1cm}
\hrulefill

\vspace{1cm}
\begin{center}
\vspace{0.2cm}
{\bf\large A discontinuous Galerkin method for optimal control of the obstacle problem}
\vspace{0.2cm}\\
{\bf\large Rohit Khandelwal$^1$}  

\vspace{0.5cm}
{\large $^1$George Mason University
}
\end{center}

\vspace{0.5cm}
This article provides a priori error estimates for an optimal control problem constrained by an elliptic obstacle problem where the finite element discretization is carried out using the symmetric interior penalty discontinuous Galerkin method. The main proofs are based on the improved $L_2$-error estimates for the obstacle problem, the discrete maximum principle, and a well known quadratic growth property. All the existing results require restrictive assumptions on mesh which is not assumed here. Quasi-optimal rate of convergence is derived for both state and control variables, in a realistic, locally distributed optimal control setting.


\vspace{1cm}
\hrulefill

\vspace{1cm}
\begin{center}
\vspace{0.2cm}
{\bf\large Mathematical modelling of the impacts of screening and vaccination in HPV}
\vspace{0.2cm}\\
{\bf\large Hyunah Lim$^1$}  

\vspace{0.5cm}
{\large $^1$University of Maryland, College Park
}
\end{center}

\vspace{0.5cm}
Cervical cancer is one of the most frequent cancers women may suffer, which caused about 342,000 death cases in 2020 worldwide. It is known that cervical cancer is highly attributable to HPV, while HPV has no medical treatment leading to cure of infection but realistic prevention only available from vaccine. However, HPV vaccines are overly expensive for individuals, and many countries including South Korea has been vaccinating young female population only for cost-effectiveness. We wish to study the impacts of different vaccine strategies such as vaccinating males and females altogether in South Korea, and also experiment on the impact of detection in early cancer stages for females through screening, to find the ideal strategy for our society to take. Although we experiment on South Korean data, our model is not regional, and we expect it to contribute to researches similar in different countries also.


\vspace{1cm}
\hrulefill

\vspace{1cm}
\begin{center}
\vspace{0.2cm}
{\bf\large Adaptive importance sampling for topology optimization in uncertain environment}
\vspace{0.2cm}\\
{\bf\large Anton Malandii$^1$}  

\vspace{0.5cm}
{\large $^1$Stony Brook University
}
\end{center}

\vspace{0.5cm}
The paper addresses the problem of topology optimization in the presence of uncertainty. The risk of failure of a structure is measured by conditional value-at- risk (CVaR). For the considered structures, estimation of CVaR requires a very large number of scenarios of failures. It is not possible to evaluate the performance of the model for every scenario because of the high computational costs associated with solving partial differential equations (PDEs). To tackle this computational complexity, we use an adaptive importance sampling approach. We efficiently estimate and minimize the CVaR by sampling failure scenarios depending on their impact. By carefully selecting critical failure scenarios, we significantly reduce the computational burden associated with solving PDE for each scenario. We validated the suggested approach with two numerical case studies.


\vspace{1cm}
\hrulefill

\vspace{1cm}
\begin{center}
\vspace{0.2cm}
{\bf\large Physics-Informed Deep Learning for Traffic State Estimation: A Survey and the Outlook
}
\vspace{0.2cm}\\
{\bf\large  Zhaobin Mo$^1$}  

\vspace{0.5cm}
{\large $^1$Columbia University
}
\end{center}

\vspace{0.5cm}
For its robust predictive power (compared to pure physics-based models) and sample-efficient training (compared to pure deep learning models), physics-informed deep learning (PIDL), a paradigm hybridizing physics-based models and deep neural networks (DNNs), has been booming in science and engineering fields. One key challenge of applying PIDL to various domains and problems lies in the design of a computational graph that integrates physics and DNNs. In other words, how the physics is encoded into DNNs and how the physics and data components are represented. We offer an overview of a variety of architecture designs of PIDL computational graphs and how these structures are customized to traffic state estimation (TSE), a central problem in building a transportation digital twin. When observation data, problem type, and goal vary, we demonstrate potential architectures of PIDL computational graphs and compare these variants using the same real-world dataset.

\vspace{1cm}
\hrulefill

\vspace{1cm}
\begin{center}
\vspace{0.2cm}
{\bf\large Mathematical modelling of the impacts of screening and vaccination in HPV}
\vspace{0.2cm}\\
{\bf\large Soyoung Park$^1$}  

\vspace{0.5cm}
{\large $^1$University of Maryland
}
\end{center}

\vspace{0.5cm}
Cervical cancer is one of the most frequent cancers women may suffer, which caused about 342,000 death cases in 2020 worldwide. It is known that cervical cancer is highly attributable to HPV, while HPV has no medical treatment leading to cure of infection but realistic prevention only available from vaccine. However, HPV vaccines are overly expensive for individuals, and many countries including South Korea has been vaccinating young female population only for cost-effectiveness. We wish to study the impacts of different vaccine strategies such as vaccinating males and females altogether in South Korea, and also experiment on the impact of detection in early cancer stages for females through screening, to find the ideal strategy for our society to take. Although we experiment on South Korean data, our model is not regional, and we expect it to contribute to researches similar in different countries also.


\vspace{1cm}
\hrulefill
\vspace{1cm}
\begin{center}
\vspace{0.2cm}
{\bf\large A digital cerebral vasculature database for stroke management}
\vspace{0.2cm}\\
{\bf\large Aseem Milind Pradhan$^1$}  

\vspace{0.5cm}
{\large $^1$George Mason University
}
\end{center}

\vspace{0.5cm}
Human cerebral collateral circulation acts as a vital safeguard against interruptions in blood flow, especially during strokes. The quality of collateral circulation significantly affects stroke treatment outcomes. However, assessing this circulation is subjective due to limitations in current metrics, and advanced imaging techniques cannot capture detailed collateral structures, leaving us with an incomplete understanding of the cerebral vasculature. 

In this study, we use a constructive constrained optimization algorithm to extend arterial trees reconstructed from MR images, generating realistic vascular models spanning from the aortic arch to cerebral arteries as small as 50 micrometers. These models can incorporate various Circle of Willis configurations and collateral vessel distributions to accommodate individual subject differences. Consequently, we aim to establish a comprehensive database of cerebral vasculatures (i.e., a virtual population) by combining multiple segmented geometries with varied extended arterial trees. We have developed a distributed lumped parameter model coupled with a mass transport solver to enable in silico angiographies for such large vascular networks. These virtual angiograms can be compared to patient-specific angiograms, identifying a model within the population (digital twin) that best approximates the patient's cerebral vasculature. 

It is essential to note that this model differs from a traditional digital twin, as information transfer between the virtual model and the patient occurs infrequently, typically only during their imaging sessions. Nevertheless, these virtual geometries facilitate blood flow simulations providing data inaccessible through conventional imaging. They can predict the impact of interventions like blockages, hemorrhages, or catheter insertion on a patient's hemodynamics. This innovative approach offers promising prospects to enhance stroke management by providing a deeper understanding of cerebral vasculature and the consequences of various interventions.

\vspace{1cm}
\hrulefill

\vspace{1cm}
\begin{center}
\vspace{0.2cm}
{\bf\large Accelerating the computation of tensor $Z$-eigenvalues}
\vspace{0.2cm}\\
{\bf\large Rhea Shroff$^1$}  

\vspace{0.5cm}
{\large $^1$University of Florida
}
\end{center}

\vspace{0.5cm}
Efficient solvers for tensor eigenvalue problems are important tools for the analysis of higher-order data sets. Here we introduce, analyze and demonstrate an extrapolation method to accelerate the widely used shifted symmetric higher order power method for tensor $Z$-eigenvalue problems. We analyze the asymptotic convergence of the method, determining the range of extrapolation parameters that induce acceleration, as well as the parameter that gives the optimal convergence rate. We then introduce an automated method to dynamically approximate the optimal parameter, and demonstrate it's efficiency when the base iteration is run with either static or adaptively set shifts. Our numerical results on both even and odd order tensors demonstrate the theory and show we achieve our theoretically predicted acceleration.

\vspace{1cm}
\hrulefill

\vspace{1cm}
\begin{center}
\vspace{0.2cm}
{\bf\large Predictive system resilience modeling}
\vspace{0.2cm}\\
{\bf\large Priscila Silva$^1$}  

\vspace{0.5cm}
{\large $^1$University of Massachusetts Dartmouth
}
\end{center}

\vspace{0.5cm}
Resilience is the ability of a system to respond, absorb, adapt, and recover from a disruptive event. Dozens of metrics to quantify resilience have been proposed in the literature. However, fewer studies have proposed models to predict these metrics or the time at which a system will be restored to its nominal performance level after experiencing degradation. This talk presents three alternative approaches to model and predict performance and resilience metrics with techniques from reliability engineering, including (i) bathtub-shaped hazard functions, (ii) mixture distributions, and (iii) a model incorporating covariates related to the intensity of events that degrade performance as well as efforts to restore performance. These models are general, and therefore, are highly applicable to digital twins in several ways, as they can enhance the understanding, performance, and management of physical systems or processes represented by digital twins, as well as critical infrastructure, cybersecurity, and machine learning systems. However, data in these areas is not shared widely. For this reason, historical data sets on job losses during seven different recessions in the United States are used to assess the predictive accuracy of these approaches, including the recession that began in 2020 due to COVID-19. Goodness of fit measures and confidence intervals as well as interval-based resilience metrics are computed to assess how well the models perform on the data sets considered. The results suggest that both bathtub-shaped functions and mixture distributions can produce accurate predictions for data sets exhibiting V, U, L, and J shaped curves, but that W and K shaped curves that respectively experience multiple shocks, deviate from the assumption of a single decrease and subsequent increase, or suffers a sudden drop in performance cannot be characterized well by either of those classes proposed. In contrast, the model incorporating covariates is capable of tracking all of types of curve noted above very well, including W and K shaped curves such as the two successive shocks the U.S. economy experienced in 1980 and the sharp degradation in 2020. Moreover, covariate models outperform the simpler models on all of the goodness of fit measures and interval-based resilience metrics computed for all seven data sets considered. These results suggest that classical reliability modeling techniques such as bathtub-shaped hazard functions and mixture distributions are suitable for modeling and prediction of some resilience curves possessing a single decrease and subsequent recovery, but that covariate models to explicitly incorporate explanatory factors and domain specific information are much more flexible and achieve higher goodness of fit and greater predictive accuracy. Thus, the covariate modeling approach provides a general framework for data collection and predictive modeling for a variety of resilience curves.


\vspace{1cm}
\hrulefill

\vspace{1cm}
\begin{center}
\vspace{0.2cm}
{\bf\large Community robustness under edge addition in synthetic and empirical temporal networks}
\vspace{0.2cm}\\
{\bf\large  Moyi Tian$^1$}  

\vspace{0.5cm}
{\large $^1$Brown University
}
\end{center}

\vspace{0.5cm}
Communities often represent key structural and functional clusters in complex networks. Many real-world networks exhibit the property of edge densification over time and observed networks may erroneously present false edges. Therefore, it is important to understand how communities change as the network changes and examine the limits of the robustness of communities under expanding the network by edge addition. We study the effect of edge addition on both real-world temporal networks and its digital twins, the simulated dynamic networks using synthetic data with controlled edge-addition strategy. We present experimental methodologies designed for synthetic and empirical cases and demonstrate results from applying four state-of-the-art community detection algorithms, i.e., Infomap, Label propagation, Leiden, and Louvain, on various network datasets. Our findings suggest that the robustness of communities depends strongly on the choice of detection method.


\vspace{1cm}
\hrulefill

\vspace{1cm}
\begin{center}
\vspace{0.2cm}
{\bf\large Large deviation theory-based adaptive importance sampling for rare events in high dimensions}
\vspace{0.2cm}\\
{\bf\large Shanyin Tong$^1$}  

\vspace{0.5cm}
{\large $^1$Columbia University
}
\end{center}

\vspace{0.5cm}
Rare and extreme events like hurricanes, energy grid blackouts, dam breaks, earthquakes, and pandemics are infrequent but have severe consequences. Because estimating the probability of such events can inform strategies that mitigate their effects, scientists must develop methods to study the distribution tail of these occurrences. However, calculating small probabilities is hard, particularly when involving complex dynamics and high-dimensional random variables. In this poster, I will discuss our proposed method for the accurate estimation of rare event or failure probabilities for expensive-to-evaluate numerical models in high dimensions. The proposed approach combines ideas from large deviation theory and adaptive importance sampling. The importance sampler uses a cross-entropy method to find an optimal Gaussian biasing distribution, and reuses all samples made throughout the process for both, the target probability estimation and for updating the biasing distributions. Large deviation theory is used to find a good initial biasing distribution through the solution of an optimization problem. Additionally, it is used to identify a low-dimensional subspace that is most informative of the rare event probability. This subspace is used for the cross-entropy method, which is known to lose efficiency in higher dimensions. The proposed method does not require smoothing of indicator functions nor does it involve numerical tuning parameters. We compare the method with a state-of-the-art cross-entropy-based importance sampling scheme using three examples: a high-dimensional failure probability estimation benchmark, a problem governed by a diffusion partial differential equation, and a tsunami problem governed by the time-dependent shallow water system in one spatial dimension.


\vspace{1cm}
\hrulefill

\vspace{1cm}
\begin{center}
\vspace{0.2cm}
{\bf\large Neural Network Approaches for Parameterized Optimal Control}
\vspace{0.2cm}\\
{\bf\large Deepanshu Verma$^1$}  

\vspace{0.5cm}
{\large $^1$Emory University
}
\end{center}

\vspace{0.5cm}
We present numerical approaches for deterministic, finite-dimensional optimal control problems whose dynamics depend on an unknown or uncertain parameter. The objective is to amortize the solution over a set of relevant parameters in an offline stage to enable rapid decision-making and be able to react to changes in the parameter in the online stage.  To tackle the curse of dimensionality arising when the state and or parameter dimension are high-dimensional, we represent the policy using neural networks. We compare two training paradigms:  First, our model-driven approach leverages the dynamics and definition of the objective function to learn the value function of the parameterized optimal control problem and obtain the policy using a feedback form.  Second, we use actor-critic reinforcement learning to approximate the policy in a data-driven way. Through a two-dimensional convection-diffusion equation, featuring high-dimensional state and parameter spaces, we investigate the accuracy, efficiency, and scalability of both training paradigms. While both paradigms lead to a reasonable approximation of the policy, the model-driven approach is more accurate and reduces the number of PDE solves significantly.

\vspace{1cm}
\hrulefill

%
%
%




\end{document}